\documentclass[12pt]{amsart}
\usepackage{enumerate}
\input{amssymb.sty}
\title{About the QWEP conjecture}
\author{Narutaka OZAWA}
\address{Department of Mathematical Sciences,
University of Tokyo, Komaba, 153-8914}
\email{narutaka@ms.u-tokyo.ac.jp}
\date{June 1, 2003}
\thanks{The author was supported by the JSPS 
Postdoctoral Fellowships for Research Abroad}

\newtheorem{thm}{Theorem}[section]
\newtheorem{prop}[thm]{Proposition}
\newtheorem{lem}[thm]{Lemma}
\newtheorem{cor}[thm]{Corollary}
\theoremstyle{definition}
\newtheorem{defn}[thm]{Definition}
\newcommand{\G}{\Gamma} 
\newcommand{\U}{{\mathcal U}} 
\newcommand{\e}{\varepsilon} 
\newcommand{\p}{\varphi}
\newcommand{\N}{{\mathbb N}} 
\newcommand{\Z}{{\mathbb Z}} 
\newcommand{\B}{{\mathbb B}} 
\newcommand{\M}{{\mathbb M}} 
\newcommand{\K}{{\mathbb K}} 
\newcommand{\R}{{\mathbb R}} 
\newcommand{\C}{{\mathbb C}} 
\newcommand{\F}{{\mathbb F}} 
\newcommand{\FF}{{\mathcal F}} 
\newcommand{\hh}{{\mathcal H}} 
\newcommand{\hhh}{\widehat{{\mathcal H}}} 
\newcommand{\id}{\mathrm{id}} 
\newcommand{\re}{\mathop{\mathrm{Re}}}
\newcommand{\im}{\mathop{\mathrm{Im}}}
\newcommand{\dist}{\mathop{\mathrm{d}}\nolimits} 
\newcommand{\tr}{\mathop{\mathrm{tr}}\nolimits} 
\newcommand{\Tr}{\mathop{\mathrm{Tr}}\nolimits} 
\newcommand{\Ad}{\mathop{\mathrm{Ad}}\nolimits} 
\newcommand{\op}{\mathrm{op}} 
\newcommand{\CB}{\mathop{\mathrm{CB}}\nolimits} 
\newcommand{\cb}{\mathrm{cb}} 
\newcommand{\ball}{\mathop{\mathrm{Ball}}\nolimits} 
\newcommand{\cball}{\mathop{\overline{\mathrm{Ball}}}\nolimits} 
\newcommand{\mult}{\mathop{\mathcal{M}}\nolimits} 
\newcommand{\ip}[1]{\langle#1\rangle} 
\begin{document}
\begin{abstract}
This is a detailed survey on the QWEP conjecture 
and Connes' embedding problem. 
Most of contents are taken from Kirchberg's paper 
[Invent.\ Math.\ \textbf{112} (1993)]. 
\end{abstract}
\maketitle
\section{Introduction}
Following Kirchberg \cite{klp}, we prove that 
several important conjectures arising in several branches 
of operator algebras are in fact equivalent. 

\begin{thm}\label{conjectures}
The following conjectures are equivalent. 
\begin{enumerate}[\rm (i)]
\item
We have 
$C^*\F_\infty\otimes_{\max}C^*\F_\infty
=C^*\F_\infty\otimes_{\min}C^*\F_\infty$. 

\item
The predual of any (separable) von Neumann algebra is 
finitely representable in the trace class $S_1$. 

\item
Any separable $\mathrm{II}_1$-factor is a subfactor of 
the ultrapower $R^\omega$ of 
the hyperfinite $\mathrm{II}_1$-factor $R$. 
\end{enumerate}
\end{thm}

Let $A$ and $B$ be $C^*$-algebras and consider
the problem to introduce a norm on their algebraic 
tensor product $A\otimes B$. 
Although there are presumably many $C^*$-norms on 
the $*$-algebra $A\otimes B$, there are two distinguished norms; 
the minimal one and the maximal one. 
(A rather surprising theorem of Takesaki \cite{takesakiI} states 
that the minimal tensor norm is indeed the smallest $C^*$-norm 
on $A\otimes B$.) 
These $C^*$-norms have nice functorial properties; 
both tensor norms tensorize completely positive maps and 
persist in passing to the second dual, 
and moreover, the minimal tensor norm is injective 
while the maximal tensor norm is projective. 
Thus, it is natural to ask whether or not the minimal (resp.\ maximal) 
tensor norm is the unique $C^*$-norm with these functoriality. 
Since injective (resp.\ projective) $C^*$-norm comes from 
a fixed $C^*$-norm on $\B(\ell_2)\otimes\B(\ell_2)$ 
(resp.\ $C^*\F_\infty\otimes C^*\F_\infty$), 
this problem is equivalent to asking 
whether there is only one $C^*$-norm on 
$\B(\ell_2)\otimes\B(\ell_2)$ 
(resp.\ $C^*\F_\infty\otimes C^*\F_\infty$). 
Probably led by this consideration (cf.\ \cite{kicm}), 
Kirchberg conjectured that the answer would be yes 
for both cases. 
However, the problem on $\B(\ell_2)\otimes\B(\ell_2)$
was solved negatively by Junge and Pisier \cite{jp}. 
The problem on $C^*\F_\infty\otimes C^*\F_\infty$ remains outstanding. 

The conjecture (ii) seems old (cf.\ \cite{hp}). 
Since several analytical structures of Banach spaces 
depend only on those of finite dimensional subspaces, 
it is desirable to know whether any noncommutative $L_1$-spaces 
are finitely representable in the trace class $S_1$. 
Many problems arising in operator spaces and noncommutative $L_p$ spaces 
are also connected to finite representability 
(as an operator space, e.g. \cite{junge}\cite{ejr}\cite{jr}). 

The problem (iii) was casually raised by Connes \cite{connes} 
and is drawing increased attention in recent years
in connection with free probability theory. 
Haagerup \cite{haagerupinv} showed 
a remarkable theorem that invariant subspaces exist 
for a large class of operators in $\mathrm{II}_1$-factors 
that are embeddable in $R^\omega$. 
Thus, a positive answer would nicely complement Haagerup's theorem. 
The problem (iii) is equivalent to that microstates always exist. 
Microstates are used to define free entropy $\chi$ introduced 
by Voiculescu \cite{voiculescuI}. 
Besides its own interest, free entropy has a number of 
important applications to theory of von Neumann algebras. 
Voiculescu \cite{voiculescuV} later introduced another 
free entropy $\chi^*$ and asked the ``unification problem'' 
whether they actually coincide. 
A negative answer to (iii) would imply this is not the case. 
See Voiculescu's survey \cite{voiculescusurvey} for details. 

The problem (iii) may well be connected to geometric group theory 
in deep and important ways. 
A positive answer would imply that all countable discrete groups 
are hyperlinear and hence would refute the famous ``theorem'' of 
Gromov claiming that a proposition which holds for all countable 
discrete groups is either trivial or false. 
Here, we say a group is hyperlinear if it is embeddable 
into the unitary group $\U(R^\omega)$ of $R^\omega$. 
The class of hyperlinear groups is closed under several natural operations 
and contains all amenable groups and all residually finite groups. 
Since many exotic groups (or ``monsters''), 
such as periodic simple groups with Kazhdan's property $\mathrm{(T)}$  
and non-uniformly embeddable groups of Gromov \cite{gromov} \cite{grorandom}, 
are arising as limits of hyperbolic groups, 
it is particularly interesting to know whether 
all hyperbolic groups are hyperlinear. 
In particular, if such simple property $\mathrm{(T)}$ groups 
or non-uniformly embeddable groups are not hyperlinear, 
then there would exist non-hyperlinear hyperbolic groups. 
Whether all hyperbolic groups are residually finite (and hence hyperlinear) 
is one of the major open problem in geometric group theory. 

It is not completely understood for which $C^*$-algebra $A$ the semigroup 
$\mathrm{Ext}(A)$ is actually a group. 
Kirchberg \cite{klp} showed the $\mathrm{Ext}$ semigroup 
of the cone (or the suspension) over $A$ is a group 
if and only if $A$ has the LLP (local lifting property). 
Since a positive answer would imply that exact and non-nuclear 
$C^*$-algebras cannot have the LLP, it would follow 
$\mathrm{Ext}(\mathrm{cone}(A))$ (and $\mathrm{Ext}(S(A))$) 
never a group when $A$ is exact and non-nuclear. 

\begin{proof}[Proof of Theorem \ref{conjectures}]
The equivalence (i)$\Leftrightarrow$(ii) follows from 
Propositions \ref{conj2}, \ref{perm1} and Corollary \ref{frep}. 
The implication (ii)$\Rightarrow$(iii) follows from Corollary \ref{emb}. 
Now, we assume (iii) and prove (i). 
The assumption (iii) is equivalent to that 
all extremal traces (and \textit{a fortiori} all traces) 
of $C^*\F_\infty$ satisfy the conditions in Theorem \ref{trace}. 
This implies that all separable von Neumann algebras 
of type $\mathrm{II}_1$ are QWEP (cf. Corollary \ref{emb}). 
We complete the proof by Proposition \ref{perm1} 
and the well-known Takesaki's theorem \cite{takesakiII} that 
any separable von Neumann algebra is $*$-isomorphic to 
a cp complemented von Neumann subalgebra in 
a separable semifinite von Neumann algebra. 
\end{proof}

\noindent\textbf{Acknowledgment.} 
The author would like to thank Professors Nate Brown 
and Gilles Pisier for valuable comments. 
This survey was written while the author was visiting 
the University of California at Berkeley under 
the support of the Japanese Society for the Promotion of Science
Postdoctoral Fellowships for Research Abroad.
\medskip

\noindent\textbf{Notations.}
Unless otherwise stated, all $C^*$-algebras in this paper 
are assumed unital (except for ideals) and maps between them 
are assumed linear. 
We will write 
$A,B,C,\ldots$ for unital $C^*$-algebras, 
$J,\ldots$ for ideals (which always mean closed two-sided ideals) 
in $C^*$-algebras, 
$M,N,\ldots$ for von Neumann algebras, 
$R$ for the hyperfinite $\mathrm{II}_1$-factor, 
$\F$ (resp.\ $\F_\infty$) for a free group 
(resp.\ on countably many generators), 
$\B(\hh)$ (resp. $\K(\hh)$) for the $C^*$-algebra of 
all bounded linear (resp. compact) operators on a Hilbert space $\hh$, 
$\M_n$ for the $n$ by $n$ full matrix algebra. 

\section{Preliminary Background}

In this section, we collect some basic results on positive maps. 
We will omit proofs if they are found in standard references. 

For $C^*$-algebras $A$ and $B$, 
we denote by $A\otimes B$ (resp.\ $A\otimes_{\min}B$, 
$A\otimes_{\max}B$) the algebraic (resp. minimal, maximal) 
tensor product. 
For von Neumann algebras $M$ and $N$, 
we denote by $M\bar{\otimes}N$ the von Neumann tensor product.
For the definitions and basic properties 
of these tensor product, 
we refer the reader to Chapter IV in Takesaki's book \cite{takesakiI} 
and Appendix (T) in Wegge-Olsen's book \cite{wegge-olsen}. 

A general reference for the positive maps 
is Sections 2 and 3 in Paulsen's book \cite{paulsen}. 
For simplicity, we only deal with 
unital $C^*$-algebras and unital maps, 
although this restriction is inessential. 

\begin{lem}\label{paulsen}
A unital map between unital $C^*$-algebras 
is positive if and only if it is contractive. 
A positive map from a commutative $C^*$-algebra 
into a (possibly noncommutative) $C^*$-algebra 
is automatically cp. 
\end{lem}

For the proof, 
see Corollary 2.9, Proposition 3.6 and Theorem 3.11 in \cite{paulsen}. 
The following two theorems are fundamental. 
The first one is due to Arveson \cite{arveson} for cp maps 
and Wittstock for cb maps. 
Recall that a unital self-adjoint subspace of a $C^*$-algebra 
is called an operator system. 
This theorem means that $\B(\hh)$ is \textit{injective} 
in the category of operator systems with ucp maps 
(resp.\ operator spaces with complete contractions). 

\begin{thm}\label{hb}
Let $X\subset B$ be an operator space. 
Then any cb map $\p\colon X\to\B(\hh)$ extends to 
a cb map $\tilde{\p}\colon B\to\B(\hh)$ with 
$\|\tilde{\p}\|_{\cb}=\|\p\|_{\cb}$. 
In particular, a ucp map from an operator system 
into $\B(\hh)$ extends to a ucp map. 
\end{thm}

The second one is the Stinespring dilation theorem. 

\begin{thm}
Let $\p\colon A\to\B(\hh)$ be a ucp map. 
Then, there is a representation $\pi\colon A\to\B(\hhh)$ 
and an isometry $V\colon\hh\to\hhh$ such that 
$\pi(A)V\hh$ is dense in $\hhh$ and that 
$\p(a)=V^*\pi(a)V$ for $a\in A$. 
Moreover, the triple $(\pi,\hhh,V)$ 
is unique up to unitary equivalence. 
\end{thm}

\begin{cor}\label{max}
If $\p\colon A\to B$ is a ucp map, 
then $\p\otimes\id_C\colon A\otimes C\to B\otimes C$
is continuous w.r.t.\ the minimal (resp.\ maximal) tensor norms. 
\end{cor}
\begin{proof}
We only prove this corollary for the maximal tensor product. 
Take a faithful representation $B\otimes_{\max}C\subset\B(\hh)$ 
and let $(\pi,\hhh,V)$ be the Stinespring triplet 
for $\p\colon A\to\B(\hh)$. 
We claim that for $x\in C$, the operator $\rho(x)$ defined by 
$$\rho(x)\colon\pi(A)V\hh\ni\sum_j\pi(a_j)V\xi_j
\longmapsto\sum_j\pi(a_j)Vx\xi_j\in\hhh$$
is well-defined and $\|\rho(x)\|\le\|x\|$. 
Indeed, putting $\xi=[\xi_1,\ldots,\xi_n]^T\in\ell_2^n\otimes\hh$, 
we have 
\begin{align*}
\|\sum_j\pi(a_j)Vx\xi_j\|^2 &=\sum_{i,j}(\p(a_i^*a_j)x\xi_j,x\xi_i)\\
 &= ((1\otimes x)^*[\p(a_i^*a_j)](1\otimes x)\xi,\xi)\\
 &\le \|x\|^2([\p(a_i^*a_j)]\xi,\xi) 
 =\|x\|^2\|\sum_j\pi(a_j)V\xi_j\|^2,
\end{align*}
where we used the fact that $[\p(a_i^*a_j)]_{i,j}\in\M_n(B)$ 
is positive and commutes with $1\otimes x$ for $x\in C\subset B'$. 
Since $\pi(A)V\hh\subset\hhh$ is dense, 
$\rho$ gives rise to a representation of 
$C$ on $\hhh$ whose range commutes with $\pi(A)$. 
Hence, $\pi\times\rho$ is a representation of $A\otimes_{\max}C$ 
such that 
$(\p\otimes\id_C)(a\otimes x)=V^*(\pi\times\rho)(a\otimes x)V$
for $a\in A$ and $x\in C$. 
This completes the proof. 
The above construction is taken from \cite{arveson}. 
\end{proof}

The Jordan product $\circ$ on a $C^*$-algebra $A$ is defined as 
$a\circ b=(ab+ba)/2$ for $a$ and $b$ in $A$. 
This equips $A$ with 
a commutative non-associative $*$-algebra structure. 
A self-adjoint subspace in a $C^*$-algebra which is closed under 
the Jordan product will be simply called a Jordan algebra. 
A positive linear map between Jordan algebras is called 
a Jordan morphism if it preserves the Jordan product. 
Jordan morphisms are necessarily contractive 
(cf.\ Lemma \ref{paulsen}). 

The following Schwarz type inequalities are due to 
Kadison \cite{kadison} for positive maps and Choi \cite{choi} for cp maps.

\begin{cor}\label{schwarz}
Let $\p\colon A\to\B(\hh)$ be a ucp (resp.\ unital positive) map.
Then, we have $\p(a^*a)\geq\p(a)^*\p(a)$ 
(resp.\ $\p(a^*\circ a)\geq\p(a)^*\circ\p(a)$)
for every $a\in A$. 
\end{cor}
\begin{proof}
First, let $\p$ be a ucp map and $(\pi,\hhh,V)$ 
be the Stinespring triplet. 
Then, we have 
$$\p(a^*a)-\p(a)^*\p(a)=V^*\pi(a)^*(1-VV^*)\pi(a)V\geq 0$$
for every $a\in A$. 
This proves the case where $\p$ is cp. 

Now, let $\p$ be a unital positive map. 
It follows from the cp case that 
$\p(b^2)\geq\p(b)^2$ for every self-adjoint element $b$ in $A$ 
since the restriction of $\p$ to the commutative $C^*$-subalgebra 
generated by $b$ is automatically cp. 
Hence, denoting by $\re a=(a^*+a)/2$ and $\im a=i(a^*-a)/2$, we have 
$$\p(a^*\circ a)=\p((\re a)^2+(\im a)^2)
\geq\p(\re a)^2+\p(\im a)^2=\p(a)^*\circ\p(a)$$
for every $a\in A$.
\end{proof}

\begin{cor}\label{mdom}
Let $\p\colon A\to\B(\hh)$ be a ucp (resp.\ unital positive) map. 
For $a\in A$, we have 
\begin{align*}
\p(a^*a)=\p(a)^*\p(a)&\Rightarrow\p(xa)=\p(x)\p(a)\mbox{ for all }x\in A\\
\p(aa^*)=\p(a)\p(a)^*&\Rightarrow\p(ax)=\p(a)\p(x)\mbox{ for all }x\in A\\
\mbox{(resp.\ }\p(a^*\circ a)=\p(a)^*\circ\p(a) &\Rightarrow 
\p(x\circ a)=\p(x)\circ\p(a)\mbox{ for all }x\in A\mbox{).}
\end{align*}
Moreover, the subspace 
$C=\{ a\in A : \p(a^*a)=\p(a)^*\p(a),\ \p(aa^*)=\p(a)\p(a)^*\}$
(resp.\ $C=\{ a\in A : \p(a^*\circ a)=\p(a)^*\circ\p(a)\}$) 
is a $C^*$-subalgebra (resp.\ a Jordan subalgebra) of $A$.  
\end{cor}

\begin{proof}
We denote by $\cdot$ the usual product (resp.\ the Jordan product) 
throughout this proof. 
Let $a\in A$ be so that $\p(a^*\cdot a)=\p(a)^*\cdot\p(a)$ 
and $x\in A$ be arbitrary. 
By Corollary \ref{schwarz}, we have for any $t\in\R$,
\begin{align*}
0 &\le\p((ta+x)^*\cdot(ta+x))-\p(ta+x)^*\cdot\p(ta+x)\\ 
 &=t\big(\p(a^*\cdot x+x^*\cdot a)-\p(a)^*\cdot\p(x)
  -\p(x)^*\cdot\p(a)\big)+\p(x^*\cdot x)-\p(x)^*\cdot\p(x).
\end{align*}
Since $t\in\R$ is arbitrary, it follows that 
$\p(a^*\cdot x+x^*\cdot a)=\p(a)^*\cdot\p(x)+\p(x)^*\cdot\p(a)$.
Replacing $x$ with $ix$, we obtain 
$\p(a^*\cdot x-x^*\cdot a)=\p(a)^*\cdot\p(x)-\p(x)^*\cdot\p(a)$.
Combining them, we obtain 
$\p(x^*\cdot a)=\p(x)^*\cdot\p(a)$
for all $x\in A$. 
This proves the first half.

We only prove that $C$ is a Jordan subalgebra for unital positive $\p$. 
It is clear that $C$ is a self-adjoint subspace of $A$. 
Hence, it is sufficient to show that $C_{\mbox{sa}}=\{a\in C : a=a^*\}$ 
is closed under the Jordan product. 
Since the usual product and the Jordan product 
coincide on a commutative subalgebra, 
we have $a^2\in C_{\mbox{sa}}$ provided that $a\in C_{\mbox{sa}}$. 
Now, the claim that $a\circ b\in C_{\mbox{sa}}$ for 
any $a,b\in C_{\mbox{sa}}$
follows from the first half of this corollary and the equation
$$(a\circ b)^2=a\circ(b\circ(a\circ b))+\frac{1}{2}(a^2\circ b)\circ b
 -\frac{1}{2}a^2\circ b^2.$$
This completes the proof. 
\end{proof}

\begin{defn}\label{mdomdef}
We say the subalgebra $C$ in Corollary \ref{mdom} is 
the \emph{multiplicative domain} for $\p$. 
\end{defn}
By Corollary \ref{schwarz}, 
$a\in A$ is in the multiplicative domain for $\p$ 
provided that that $\|a\|=1$ and $\p(a)$ is a unitary. 
Therefore, if $\p\colon A\to B$ is a ucp (resp.\ unital positive) map 
which maps the closed unit ball of $A$ 
onto the closed unit ball of $B$, 
then $\p$ is surjective on the multiplicative domain for $\p$. 

The following important corollary will be used frequently 
without mention. 

\begin{cor}\label{bimod}
Let $A\subset B$ and 
let $\p\colon B\to\B(\hh)$ be a ucp map such that 
the restriction of $\p$ to $A$ is a $*$-homomorphism $\pi$, 
then $\p$ is an $A$-bimodule map, i.e., 
$\p(axb)=\pi(a)\p(x)\pi(b)$ for $a,b\in A$ and $x\in B$. 
\end{cor}

A map $\p$ on $B$ is called a projection if $\p^2=\p$.
The following theorem is due to Tomiyama \cite{tomiyama}. 
A simple proof is found in Section 9 in Str{\u a}til{\u a}'s book \cite{stratila}
and Chapter IX in Takesaki's book \cite{takesakiII}.

\begin{thm}\label{tomiyama}
For $A\subset B$ and a projection $\p$ from $B$ onto $A$, 
the following are equivalent.

\begin{enumerate}[\rm (i)]
\item
The map $\p$ is a conditional expectation, i.e., 
the map $\p$ is an $A$-bimodule map; 
$\p(axb)=a\p(x)b$ for $a,b\in A$ and $x\in B$. 

\item
The map $\p$ is cp.

\item
The map $\p$ is contractive.
\end{enumerate}
\end{thm}

We need the structure theorem for Jordan morphisms 
due to St{\o}rmer \cite{hs}. 

\begin{thm}\label{jordan}
Let $A$ be a $C^*$-algebra and 
$\theta\colon A\to\B(\hh)$ be a Jordan morphism. 
Then, there is a central projection $p$ in 
the von Neumann algebra $M=\theta(A)''$ 
generated by $\theta(A)$ such that 
the map $A\ni a\mapsto\theta(a)p\in Mp$ 
(resp.\ $A\ni a\mapsto\theta(a)(1-p)\in M(1-p)$)
is a $*$-homomorphism (resp.\ $*$-antihomomorphism). 
\end{thm}

We review some operator space theory. 
We refer the reader to 
the books of Effros and Ruan \cite{er} and of Pisier \cite{pintro}.
Many important results in operator space theory are related to 
the dual operator space structure introduced by 
Blecher-Paulsen and Effros-Ruan 
\cite{blecher}, \cite{bp}, \cite{ernew}. 

Let $X\subset\B(\hh)$ be an operator space and 
let $X^*$ be its dual Banach space. 
For $x=[x_{ij}]_{i,j}\in\M_n(X)$, we define 
$\theta_x\colon X^*\ni f\longmapsto [f(x_{ij})]_{i,j}\in\M_n$. 
Let $\Omega$ be the union of the closed unit balls of $\M_n(X)$ 
($n=1,2,\ldots$). 
We introduce an operator space structure on $X^*$ by 
the isometric inclusion 
$$\Theta\colon X^*\ni f\longmapsto(\theta_x(f))_{x\in\Omega}
\in\prod_{x\in\Omega}\M_{n(x)}.$$
Unless otherwise stated, we always assume that the dual space $X^*$ 
is equipped with this operator space structure. 
(We note that $\M_n(X^{**})=\M_n(X)^{**}$ always holds isometrically, 
but $\M_n(X^*)=\M_n(X)^*$ is \emph{not} isometric.)
It is not hard to see 
$X^*\otimes_{\min}\M_n=\M_n(X^*)=\CB(X,\M_n)$
isometrically. 
It follows that for operator spaces $X$ and $Y$, 
there is a canonical isometric inclusion defined by 
$$X^*\otimes_{\min}Y\ni\sum_kf_k\otimes y_k\longmapsto
\sum_k f_k(\,\cdot\,)y_k\in\CB(X,Y).$$
The predual $M_*$ of a von Neumann algebra $M$ is 
equipped with the operator space structure induced 
from that of $M^*$. 
It is proved in \cite{blecher} that $M=(M_*)^*$ 
completely isometrically as one should expect. 
\section{WEP and LLP}

\begin{defn}\label{defwep}
We say a $C^*$-subalgebra $A$ in $B$ is cp complemented 
(resp.\ weakly cp complemented) in $B$ 
if there is a a ucp map 
$\p\colon B\to A$ (resp.\ $\p\colon B\to A^{**}$) 
such that $\p|_A=\id_A$. 

We say a $C^*$-algebra $B$ has the WEP (weak expectation property), 
if it is weakly cp complemented in $\B(\hh)$ for 
a faithful representation $B\subset\B(\hh)$. 
\end{defn}

Since $\B(\hh)$ is injective, the definition of the WEP 
does not depend on choice of faithful representation of $B$. 
We say a $C^*$-algebra is QWEP if it is a quotient of 
a $C^*$-algebra with the WEP. 
The QWEP conjecture states that all $C^*$-algebras are QWEP. 

If $M\subset N$ are von Neumann algebras 
with a faithful normal trace $\tau$ on $N$, 
then there is a unique trace preserving conditional 
expectation $\Phi$ from $N$ onto $M$ defined by the relation
$\tau(a\Phi(x))=\tau(ax)$ for $a\in M$ and $x\in M$. 
In particular, $M$ is cp complemented in $N$. 
Indeed, this follows from the fact that 
linear functionals on $M$ of the form $\tau(a\,\cdot\,)$ 
with $a\in M$ are dense in $M_*$.
Let $\Delta$ be a subgroup of $\G$. 
Then, the $C^*$-subalgebra in the full $C^*$-algebra $C^*\G$ 
generated by $\Delta$ is naturally $*$-isomorphic to $C^*\Delta$. 
Indeed, this follows from existence of induction. 
Moreover, $C^*\Delta$ is cp complemented in $C^*\G$ 
by the conditional expectation $\p$ defined by 
$\p(s)=0$ for $s\in\G\setminus\Delta$. 

\begin{lem}\label{compl}
For $C^*$-algebras $A\subset B$, the following are equivalent. 
\begin{enumerate}[\rm (i)]
\item 
The $C^*$-algebra $A$ is weakly cp complemented in $B$. 
\item
The second dual $A^{**}$ is cp complemented in $B^{**}$. 
\item
For any finite dimensional subspace $E\subset B$ and any $\e>0$, 
there exists a map $\p\colon E\to A$ 
such that $\|\p\|\le1+\e$ and $\p|_{E\cap A}=\id_{E\cap A}$. 
\end{enumerate}
In particular, if $A$ is weakly cp complemented 
in $B$ and $B$ has the WEP, then so does $A$. 
If $(A_i)_{i\in I}$ is a family of $C^*$-algebras with the WEP, 
then $\prod A_i$ has the WEP. 
\end{lem}
\begin{proof}
The equivalence (i)$\Leftrightarrow$(ii) is a consequence 
of the fact that any bounded linear map $\p$ from 
a Banach space $X$ into a dual Banach space $Y=(Y_*)^*$ 
uniquely extends to a weak$^*$-continuous map 
$\tilde{\p}$ from $X^{**}$ into $Y$. 
The implication (i)$\Rightarrow$(iii) follows from 
the principle of local reflexivity for Banach spaces. 
To prove (iii)$\Rightarrow$(ii), 
let $I$ be the set of pairs $(E,\e)$, 
where $E\subset B$ is finite dimensional and $\e>0$. 
The set $I$ is directed by the order relation 
$(E_1,\e_1)\le(E_2,\e_2)$ if and only if 
$E_1\subset E_2$ and $\e_1\geq\e_2$. 
By the condition (iii), for each $i=(E,\e)\in I$, 
there exists a map $\p_i\colon B\to A^{**}$ 
such that $\|\p_i|_E\|\le1+\e$ and $\p_i|_{E\cap A}=\id_{E\cap A}$. 
Then, any cluster point $\p\colon B\to A^{**}$ of 
the net $\{\p_i\}_{i\in I}$, in the point-weak$^*$ topology, 
satisfies that $\|\p\|\le 1$ and $\p|_A=\id_A$. 
It follows from Theorem~\ref{tomiyama} that the weak$^*$-continuous 
extension $\bar{\p}$ of $\p$ on $B^{**}$ is a conditional expectation 
from $B^{**}$ onto $A^{**}$. 

Finally, we observe that the condition (iii) 
is stable under direct product. 
\end{proof}

The notion of WEP was introduced by Lance \cite{lwep} where 
he showed the following theorem, which in particular shows 
nuclear $C^*$-algebras have the WEP. 

\begin{thm}\label{lance}
Let $A\subset B$. 
Then, $A$ is weakly cp complemented in $B$ if and only if
$A\otimes_{\max}C\subset B\otimes_{\max}C$ 
isometrically for any $C^*$-algebra $C$
(resp.\ for $C=C^*\F_\infty$ or $C=A^{\op}$). 
\end{thm}
\begin{proof}
We note that the canonical map 
$A\otimes_{\max}C\subset A^{**}\otimes_{\max}C$ 
is always isometric and the canonical map 
$\iota\otimes\id_C\colon A\otimes_{\max}C\to B\otimes_{\max}C$ 
is always contractive. 
Hence, the ``only if '' part follows from Corollary \ref{max}. 
To prove ``if'' part, let $A\subset\B(\hh)$ be 
the universal representation, i.e., $A''=A^{**}$, 
and let $C$ be a $C^*$-algebra which has a representation 
$\pi\colon C\to\B(\hh)$ with $\pi(C)''=A'$. 
By the assumption, the representation $A\otimes_{\max}C$ 
on $\B(\hh)$ extends to a ucp map 
$\Phi\colon B\otimes_{\max}C\to\B(\hh)$. 
Then, the map $\p\colon B\ni b\mapsto\Phi(b\otimes1)\in\B(\hh)$ is 
a ucp extension of $\id_A$. 
Since $\Phi$ is a $C$-bimodule map (cf.\ Corollary \ref{bimod}), 
we have $\p(b)\pi(x)=\Phi(b\otimes x)=\pi(x)\p(b)$ for 
$b\in B$ and $x\in C$, i.e., $\p(b)\in\pi(C)'=A^{**}$. 
This completes the proof. 
\end{proof}

The following proposition reduces many problems 
to that for separable $C^*$-algebras.

\begin{prop}\label{sep}
Let $B$ be a (non-separable) $C^*$-algebra and let 
$X\subset B$ be a separable subspace. 
Then, there is a separable $C^*$-subalgbera $A$ 
which contains $X$ and is weakly cp complemented in $B$. 
\end{prop}
\begin{proof}
Let $C=C^*\F_\infty$ and 
let $(B_i)_{i\in I}$ be the directed set of all 
separable $C^*$-subalgebras in $B$ ordered by inclusions. 
We claim that 
$\| x\|_{B\otimes_{\max}C}=\lim\| x\|_{B_i\otimes_{\max}C}$
for every $x\in B\otimes C$. 
Note that the limit in R.H.S.\ always exists since the net is decreasing 
and that the inequality $\le$ is clear. 
To prove the converse inequality, suppose that 
$\| x\|_{B_i\otimes_{\max}C}\geq 1$ for all $i\in I$.
It follows that for every $i\in I$, there are commuting representations 
$\pi_i\colon B_i\to\B(\hh_i)$ and $\rho_i\colon C\to\B(\hh_i)$ 
such that $\|(\pi_i\cdot\rho_i)(x)\|_{\B(\hh_i)}\geq1$. 
Let $\omega$ be a nontrivial ultrafilter on the directed set $I$ and 
let $\hh_\omega$ be the ultraproduct Hilbert space of $(\hh_i)_{i\in I}$. 
Then, the nets $(\pi_i)_{i\in I}$ and $(\rho_i)_{i\in I}$ give rise to 
commuting representations $\pi_\omega\colon B\to\B(\hh_\omega)$ 
and $\rho_\omega\colon C\to\B(\hh_\omega)$. 
Therefore, we have
$\| x\|_{B\otimes_{\max}C}
\geq\|(\pi_\omega\cdot\rho_\omega)(x)\|_{\B(\hh_\omega)}\geq 1,$ 
which proves the claim.

Let $A_0\subset B$ be a separable $C^*$-subalgebra containing $X$. 
Using the first part of this proof recursively, 
we can find an increasing sequence $(A_n)_{n=0}^\infty$ 
of separable $C^*$-subalgebras in $B$ such that 
$\|x\|_{B\otimes_{\max}C}=\|x\|_{A_n\otimes_{\max}C}$ 
for every $n\in\N$ and $x\in A_{n-1}\otimes C$.
Let $A$ be the closure of $\bigcup_{n=1}^\infty A_n$. 
Then, $A$ is a separable $C^*$-subalgebra such that 
$A\otimes_{\max}C\subset B\otimes_{\max}C$ isometrically. 
Hence, by Theorem \ref{lance}, $A$ is weakly cp complemented in $B$. 
\end{proof}

\begin{defn}
Let $\p\colon A\to B/J$ be a ucp map. 
We say $\p$ is ucp liftable if there is a 
ucp lifting $\psi\colon A\to B$, i.e., 
there is a ucp map $\psi$ such that $\p=\pi\psi$ 
for the quotient $\pi$ from $B$ onto $B/J$. 
We say $\p$ is locally ucp liftable if 
for any finite dimensional operator system $E$ in $A$, 
there is a ucp lifting $\psi\colon E\to B$. 

We say $A$ has the LP (lifting property) 
(resp.\ the LLP (local lifting property))
if any ucp map from $A$ into any quotient $C^*$-algebra $B/J$ 
is (resp.\ locally) ucp liftable. 
\end{defn}

The following Effros-Haagerup lifting theorem \cite{eh} 
characterizes those ucp maps which are locally ucp liftable. 
\begin{thm}\label{eh}
A ucp map $\p\colon A\to B/J$ is locally ucp liftable if and only if 
$$\p\otimes\id\colon A\otimes\B(\ell_2)
\to(B\otimes_{\min}\B(\ell_2))/(J\otimes_{\min}\B(\ell_2))$$
is continuous w.r.t.\ the minimal tensor norm on $A\otimes\B(\ell_2)$. 
\end{thm}

Existence of completely contractive local lifting immediately 
follows from operator space duality explained in Section 2, 
but making the lifting completely positive requires a technical lemma. 
See \cite{eh} or \cite{wassermann} for the proof of the following lemma. 

\begin{lem}\label{ehlem}
Let $\p\colon E\to B/J$ be a ucp map 
from a finite dimensional operator system into a quotient $C^*$-algebra. 
Suppose that for any $\e>0$, there is a lifting $\psi\colon E\to B$ 
with $\|\psi\|_{\cb}\le1+\e$.
Then, $\p$ is ucp liftable. 
\end{lem}

The $C^*$-algebra $\B(\ell_2)$ is universal in the sense that 
it contains all separable $C^*$-algebras and is injective. 
The other universal thing is the full $C^*$-algebra $C^*\F_\infty$ 
of the free group $\F_\infty$ on countably many generators since 
any separable $C^*$-algebra is a quotient of it and it has the LP. 

\begin{thm}\label{lp}
The full $C^*$-algebra $C^*\F$ of a countable free group $\F$ 
has the LP. 
\end{thm}

\begin{proof}
This proof is taken from \cite{kcom}.
We first show that a $*$-homomorphism $\theta\colon C^*\F\to B/J$ 
is ucp liftable. 
To this end, let $x_1,x_2,\ldots\in B$ be a contractive liftings 
of $\theta(s_1),\theta(s_2),\ldots$,
where $s_1,s_2,\ldots$ are the free generators of $\F$.
Then, each $x_n$ dilates to a unitary 
$$\hat{x}_n=\left[\begin{array}{cc} 
 x_n & (1-x_nx_n^*)^{1/2} \\ (1-x_n^*x_n)^{1/2} & -x_n^*
\end{array}\right] \in \M_2(B).$$ 
By universality, there is a $*$-homomorphism 
$\rho\colon C^*\F\to \M_2(B)$ with $\rho(s_n)=\hat{x}_n$. 
It is not hard to see that the $(1,1)$-corner of $\rho$ 
is a desired lifting of $\theta$. 

Now, let $\p\colon C^*\F\to B/J$ be a ucp map. 
Since $\F$ is countable, we may assume $B/J$ is separable. 
By the Kasparov-Stinespring dilation theorem \cite{lhilb}, 
there is a $*$-homomorphism 
$\theta\colon C^*\F\to\mult(\K\otimes_{\min}B/J)$ 
such that $\p(a)=\theta(a)_{11}$ for $a\in C^*\F$, 
where $x_{11}$ is the $(1,1)$-entry of $x\in\mult(\K\otimes_{\min}B/J)$. 
By the noncommutative Tietze extension theorem, 
the surjective $*$-homomorphism $\pi$ from $\K\otimes_{\min}B$ 
onto $\K\otimes_{\min}B/J$ extends to a surjective 
$*$-homomorphism $\tilde{\pi}$ between their multiplier algebras. 
Hence, by the first part of this proof, 
there is a ucp map
$\rho\colon C^*\F\to\mult(\K\otimes_{\min}B)$
such that $\theta=\tilde{\pi}\rho$. 
The ucp map $\psi\colon C^*\F\ni a\mapsto\rho(a)_{11}\in B$ 
is a desired lifting of $\p$. 
\end{proof}

The above proof seems unreasonably involved. 
Hence it would be interesting to find another proof. 
(The LLP, rather than the LP, is sufficient in most of places and 
an independent proof of the LLP for $C^*\F_\infty$ will be given 
in Theorem \ref{fb}.)
We note that a map $\p$ from a discrete group $\G$ into a $C^*$-algebra $A$ 
extends to a ucp map on the full $C^*$-algbera $C^*\G$  
if and only if it is unital and positive definite, i.e., 
$[\p(s_i^{-1}s_j)]_{i,j}$ is positive in $\M_n(A)$ 
for any $n$ and any $s_1,\ldots,s_n\in\G$. 

Since countable noncommutative free groups 
are isomorphic to subgroups of each other, 
their distinction in our story is very minor. 
The full $C^*$-algebras of free groups 
(on uncountable generators) do have the LLP. 
It seems they do not have the LP, but 
it is not known (even for the free group of continuous cardinality). 
There is a (non-separable) $C^*$-algebra which has the LLP 
but not the LP. 
Indeed, the commutative $C^*$-algbera $\ell_\infty/c_0$ is 
such an example as there is no bounded linear lifting from 
$\ell_\infty/c_0$ into $\ell_\infty$. 
(Since $\ell_\infty/c_0$ contains a $C^*$-subalgebra 
$*$-isomorphic to $c_0(I)$ with an uncountable index set $I$, 
it cannot be embedded (as a Banach space) into $\ell_\infty$ 
which has a separable predual.) 
We will see that the LP and the LLP are equivalent 
for separable $C^*$-algebras provided that 
the QWEP conjecture is true. 

\begin{cor}\label{liftid}
Let $A$ be a separable $C^*$-algebra and let 
$J\triangleleft C^*\F_\infty$ be such that $A=C^*\F_\infty/J$. 
Then, $A$ has the LP (resp. the LLP) if 
$\id_A$ 
is (resp.\ locally) ucp liftable. 
\end{cor}

There is a criteria when a locally ucp liftable map has a global 
lifting. 
This is due to Arveson \cite{anote} 
and its variant is due to Effros and Haagerup \cite{eh}.
The condition (ii) in this lemma is satisfied whenever $J$ is nuclear. 
The separability condition is essential as the quotient 
from $\ell_\infty$ onto $\ell_\infty/c_0$ is not grobally liftable. 

\begin{lem}\label{arv}
Let $p\colon A\to B/J$ be a ucp map. 
Then, $\p$ is ucp liftable provided that $A$ is separable 
and either $\mathrm{(i)}$ or $\mathrm{(ii)}$ holds. 

\begin{enumerate}[\rm (i)]
\item
The ucp map $\p$ can be approximated by 
ucp liftable maps in the point-norm topology. 

\item
The ucp map $\p$ is locally ucp liftable, and moreover 
for any finite dimensional operator systems $E\subset F\subset A$, 
any $\e>0$ and for any cp map $\theta\colon E\to J$, 
there is a cp map $\tilde{\theta}\colon F\to J$ 
with $\|\tilde{\theta}|_E-\theta\|<\e$. 
\end{enumerate}
\end{lem}

\begin{proof}
Let $\pi$ be the quotient map from $B$ onto $B/J$ 
and let $E_1\subset E_2\subset\cdots\subset A$ 
be an increasing sequence of finite dimensional operator system 
with dense union. 
We first deal with the case (ii). 
Let $\psi_n\colon E_n\to B$ be ucp maps such that 
$\|\p(a)-\pi\psi_n(a)\|<2^{-n}\|a\|$ for $a\in E_n$. 
We claim that there is a sequence of ucp maps 
$\tilde{\psi}_n\colon E_n\to B$ 
such that $\pi\tilde{\psi}_n=\pi\psi_n$ and 
$\|\tilde{\psi}_{n+1}(a)-\tilde{\psi}_n(a)\|<2^{-n+2}\|a\|$ 
for $n\in\N$ and $a\in E_n$. 

We proceed by induction. 
Suppose that $\tilde{\psi}_n$ is already constructed. 
Take $\e>0$ sufficiently small. 
Since $\|\pi(\psi_{n+1}(a)-\tilde{\psi}_n(a))\|<2^{-n+1}\|a\|$ 
for $a\in E_n$, we can find a quasicentral approximate unit 
$e\in J$ with $0\le e\le 1$ such that 
$$\|(1-e)^{1/2}(\psi_{n+1}(a)-\tilde{\psi}_n(a))(1-e)^{1/2}\|
<2^{-n+1}\|a\|$$
and moreover 
$$\|\tilde{\psi}_n(a)-\big(e^{1/2}\tilde{\psi}_n(a)e^{1/2}
+(1-e)^{1/2}\tilde{\psi}_n(a)(1-e)^{1/2}\big)\|<\e\|a\|$$
for $a\in E_n$ (cf.\ \cite{anote}). 
By the condition (ii), there is a cp map 
$\tilde{\theta}\colon E_{n+1}\to J$ such that 
$\|\tilde{\theta}(a)-e^{1/2}\tilde{\psi}_n(a)e^{1/2}\|<\e\|a\|$. 
We define $\tilde{\psi}_{n+1}$ by 
$$\tilde{\psi}_{n+1}(a)=b^{-1/2}\big(\tilde{\theta}(a)
+(1-e)^{1/2}\psi_{n+1}(a)(1-e)^{1/2}\big)b^{-1/2},$$ 
where $b=\tilde{\theta}(1)+1-e\in B$. 
Since $\|b-1\|<\e$, we have 
$$\|\tilde{\psi}_{n+1}(a)-\tilde{\psi}_n(a)\|
<(f(\e)+2^{-n+1})\|a\|$$
for $a\in E_n$, 
where $f$ is a continuous function with $f(0)=0$. 
This proves the claim. 

The desired ucp lifting $\psi\colon A\to B$ is now obtained 
by letting (for $a\in\bigcup E_n$)
$$\psi(a)=\lim_{n\to\infty}\tilde{\psi}_n(a).$$ 

The proof for the case (i) is similar but easier than the case (ii). 
Indeed, $\psi_n$ and $\tilde{\psi}_n$ in the above proof 
are defined globally on $A$ 
(hence no need to choose $\tilde{\theta}$). 
\end{proof}

The following is the Choi-Effros lifting theorem \cite{ce}. 
\begin{cor}
Let $\p\colon A\to B/J$ be a ucp map. 
Then, $\p$ is ucp liftable provided that 
$A$ is separable and $\p$ is nuclear, i.e., 
there is a sequence of ucp maps $\beta_n\colon A\to\M_{k(n)}$ 
and $\alpha_n\colon\M_{k(n)}\to B/J$ such that 
$\alpha_n\beta_n$ converges to $\p$ pointwise. 

In particular, a separable nuclear $C^*$-algebra has 
the LP. 
\end{cor}
\begin{proof}
It is well-known and not too hard to show that 
a map $\alpha$ from $\M_k$ into a $C^*$-algebra $C$ is 
cp if and only if $[\alpha(e_{ij})]_{ij}\in\M_k(C)$ is positive. 
It follows that any ucp map $\alpha\colon\M_k\to B/J$ is ucp liftable. 
The conclusion follows from this fact and Lemma \ref{arv}. 
\end{proof}

The following is due to Kirchberg \cite{klp}.

\begin{cor}\label{llplp}
Let $\p\colon A\to B/J$ be a ucp map. 
If $A$ is a separable $C^*$-algebra with the LLP and 
$B$ is QWEP, then $\p$ is ucp liftable. 
\end{cor}
\begin{proof}
We may assume $B$ has the WEP. 
To verify the condition (ii) in Lemma \ref{arv}, 
we give ourselves finite dimensional operator systems 
$E\subset F\subset A$, $\e>0$ and a cp map $\theta\colon E\to J$. 
Since $B$ has the WEP, $\theta$ extends to a ucp map 
from $A$ into $B^{**}$ which is still denoted by $\theta$. 
For a directed set $I$, let 
$$B_I=\{ (b_i)_{i\in I}\in\prod_{i\in I}B : 
\mbox{strong$^*$-}\lim_{i\in I} b_i\mbox{ exists in }B^{**}\}$$
and let $\pi\colon B_I\to B^{**}$ be the map which takes 
$(b_i)_{i\in I}$ to its limit. 
Since the adjoint-operation and product is jointly 
strong$^*$-continuous on bounded sets, 
$B_I$ is a $C^*$-algebra and $\pi$ is a $*$-homomorphism. 
Choosing the directed set $I$ appropriately, 
we may assume that $\pi$ is surjective. 
Since $A$ has the LLP, there is a ucp map $\psi\colon F\to B_I$ 
such that $\theta|_F=\pi\psi$. 
It follows that there is a net of ucp maps $\psi_i\colon F\to B$ 
such that the net $\psi_i|_E$ converges to $\theta$ in the 
pointwise weak topology. 
Taking convex combinations and multiplying by an approximate unit, 
we find a desired cp map $\tilde{\theta}\colon F\to J$. 
\end{proof}

The converse is also true. 
For the proof, see \cite{olp}.
\begin{prop}\label{liftcalkin}
A separable $C^*$-algebra $A$ has the LLP 
if and only if any ucp map from $A$ into 
the Calkin algebra $\B(\ell_2)/\K(\ell_2)$ is ucp liftable. 
\end{prop}

It might be the case that any ucp map from 
a separable $C^*$-algebra into the Calkin algebra 
$\B(\ell_2)/\K(\ell_2)$ dilates to a $*$-homomorphism 
into $\M_2(\B(\ell_2)/\K(\ell_2))$. 
If it is the case, then it would imply that 
a separable $C^*$-algebra has the LLP if and only if 
$\mathrm{Ext}(A)$ is a group. 
Recall that the (unitized) cone of a $C^*$-algebra $A$ is 
defined to be 
$\mathrm{cone}(A)=\{ f\in C([0,1],A) : f(0)\in\C1\}$. 
By an ingeneous argument on the cone, 
Kirchberg \cite{klp} proved that a separable $C^*$-algebra $A$ 
has the LLP if and only if $\mathrm{Ext}(\mathrm{cone}(A))$ 
is a group. 

Kirchberg \cite{kcom}\cite{klp} proved the following 
important theorem and its corollary. 
\begin{thm}\label{fb}
We have 
$C^*\F_\infty\otimes_{\min}\B(\ell_2)
=C^*\F_\infty\otimes_{\max}\B(\ell_2)$. 
\end{thm}

We will present an elegant and simple proof of Pisier \cite{psimple}. 
To this end, let $E_n$ be the $n$-dimensional operator space in $C^*\F_\infty$ 
spanned by $1=U_0,U_1,\ldots,U_{n-1}$, where $U_1,U_2,\ldots$ are 
the canonical generators of $C^*\F_\infty$. 
It is rather obvious that $E_n$ is canonically isometric to $\ell_1^n$, 
or equivalently, 
$\|\sum_{k=0}^{n-1}\alpha_kU_k\|=\sum_{k=0}^{n-1}|\alpha_k|$ 
for every $(\alpha_k)_{k=0}^{n-1}\in\C$.
This gives rise to the canonical one-to-one correspondence 
between an element 
$z=\sum_{k=0}^{n-1}U_k\otimes x_k\in E_n\otimes\B(\ell_2)$ 
and a map
$$\tilde{z}\colon\ell_\infty^n\ni(\alpha_k)_{k=0}^{n-1}
\longmapsto\sum_{k=0}^{n-1}\alpha_kx_k\in\B(\ell_2).$$

\begin{lem}\label{ell1}
The above operator space $E_n$ is canonically completely 
isometrically isomorphic to 
the dual operator space $\ell_1^n=(\ell_\infty^n)^*$,
or equivalently, $\| z\|_{\min}=\|\tilde{z}\|_{\cb}$
for every $z\in E_n\otimes\B(\ell_2)$.
\end{lem}
\begin{proof}
Since $(U_k)_{k=0}^{n-1}\in E_n\otimes_{\min}\ell_\infty$ 
is contractive and 
$z=(\id_{E_n}\otimes\tilde{z})((U_k)_{k=0}^{n-1})$, 
we have $\| z\|_{\min}\le\|\tilde{z}\|_{\cb}$. 
To prove the converse inequality, we give ourselves 
contractions $a_0,\ldots,a_{n-1}\in\B(\hh)$. 
Let $\hat{a}_k\in\M_2(\B(\hh))$ be their unitary dilations 
(cf.\ the proof of Theorem \ref{lp}). 
It follows that the map $\p\colon E_n\to\M_2(\B(\hh))$ 
defined by $\p(U_k)=\hat{a}_0^{-1}\hat{a}_k$, $k=0,\ldots,n-1$ 
is completely contractive since it extends to a $*$-homomorphism 
on $C^*\F_\infty$. 
Hence, the map $\theta\colon E_n\to\B(\hh)$ defined by 
$\theta(U_k)=(\hat{a}_0\p(U_k))_{11}=a_k$, $k=0,\ldots,n-1$ 
is also completely contractive.
Therefore, we have 
$\|(\id_{\B(\hh)}\otimes\tilde{z})((a_k)_{k=0}^{n-1})\|_{\min}
 =\|(\theta\otimes\id_{\B(\ell_2)})(z)\|_{\min}\le\|z\|_{\min}.$ 
Since the contraction 
$(a_k)_{k=0}^{n-1}\in\B(\hh)\otimes_{\min}\ell_\infty^n$ 
was arbitrary, we have $\|\tilde{z}\|_{\cb}\le\|z\|_{\min}$. 
\end{proof}

The dual operator space structure $\ell_1^n=(\ell_\infty^n)^*$ 
coincides with so called the maximal operator space structure \cite{pmax}. 

\begin{lem}\label{lin}
Let $X_i\subset\B(\hh_i)$ ($i=1,2$) be unital operator spaces 
and let $\p\colon X_1\to X_2$ be a unital complete isometry. 
Suppose that $X_2$ is spanned by unitaries in $\B(\hh_2)$. 
Then, $\p$ uniquely extends to a $*$-homomorphism between 
$C^*$-subalgebras $C^*(X_i)$ in $\B(\hh_i)$. 
\end{lem}
\begin{proof}
By Theorem \ref{hb}, $\p$ extends 
to a complete contraction from $\B(\hh_1)$ into $\B(\hh_2)$,
which is still denoted by $\p$.
Since $\p$ is unital, it has to be a ucp map.
Since $\p|_{X_1}$ is isometric and $X_2$ is spanned by unitaries, 
$X_1$ is contained in the multiplicative domain of $\p$ 
(cf.\ Definition \ref{mdomdef} and the following remarks). 
\end{proof}

\begin{proof}[Proof of Theorem \ref{fb}]
Thanks to Lemma \ref{lin}, it suffices to show 
that the formal identity map from $E_n\otimes_{\min}\B(\ell_2)$ 
into $C^*\F_\infty\otimes_{\max}\B(\ell_2)$ is completely contractive 
for every $n$. 
We give ourselves $z=\sum_{k=0}^{n-1}U_k\otimes x_k\in E_n\otimes\B(\ell_2)$ 
with $\|z\|_{\min}=1$. 
By Lemma \ref{ell1}, the corresponding map 
$\tilde{z}\colon\ell_\infty^n\to\B(\ell_2)$ is completely contractive. 
Hence, by the Stinespring-type theorem for cb maps, there are 
a Hilbert space $\hh$, a representation $\pi\colon\ell_\infty^n\to\B(\hh)$ 
and contractions $V,W\in\B(\ell_2,\hh)$ such that 
$\tilde{z}(f)=V^*\pi(f)W$ for $f\in\ell_\infty^n$. 
We may assume that $\hh=\ell_2$. 
Then, $a_k=\pi(\delta_k)V$ and $b_k=\pi(\delta_k)W$ are 
in $\B(\ell_2)$, and satisfy 
$x_k=a_k^*b_k$ for $k=0,\ldots,n-1$ and 
$\sum_{k=0}^{n-1}a_k^*a_k\le1$, $\sum_{k=0}^{n-1}b_k^*b_k\le1.$ 
It follows that 
\begin{align*}
\|\sum_{k=0}^{n-1} U_k\otimes & x_k\|_{C^*\F_\infty\otimes_{\max}\B(\ell_2)}
  =\|\sum_{k=0}^{n-1}(1\otimes a_k)^*(U_k\otimes b_k)
  \|_{C^*\F_\infty\otimes_{\max}\B(\ell_2)}\\
&\le \|\sum_{k=0}^{n-1}(1\otimes a_k)^*(1\otimes a_k)\|_{\max}^{1/2}
  \|\sum_{k=0}^{n-1}(U_k\otimes b_k)^*(U_k\otimes b_k)\|_{\max}^{1/2}\le 1.
\end{align*}
This shows that the formal identity from $E_n\otimes_{\min}\B(\ell_2)$ 
into $C^*\F_\infty\otimes_{\max}\B(\ell_2)$ is contractive. 
Since $\B(\ell_2)$ is stable, it is also completely contractive.
\end{proof}

We remark that the fact that $C^*\F_\infty$ has the LLP 
was not used in the proof of Theorem \ref{fb}. 
Combining this theorem with Theorems \ref{lance}, \ref{eh} 
and Corollary \ref{liftid}, we obtain the following. 

\begin{cor}\label{wepllp}
For $A$ and $B$, we have the following. 

\begin{enumerate}[\rm (i)]
\item
$A\otimes_{\min}B=A\otimes_{\max}B$ 
if $A$ has the LLP and $B$ has the WEP. 

\item
$A\otimes_{\min}\B(\ell_2)=A\otimes_{\max}\B(\ell_2)$ 
if and only if $A$ has the LLP.

\item
$C^*\F_\infty\otimes_{\min}B=C^*\F_\infty\otimes_{\max}B$ 
if and only if $B$ has the WEP. 
\end{enumerate}
\end{cor}

We remark that both the WEP and the LLP 
(and the LP for separable $C^*$-algebras) 
are stable under taking a tensor product with a nuclear $C^*$-algebra, 
and taking a crossed product by an amenable group. 
This fact easily follows from the above corollary 
(or the completely positive approximation property). 

\begin{cor}\label{lrp}
Let $A_i$ ($i=1,2$) be $C^*$-algebras 
and let $\pi_i\colon C^*\F\to A_i$ be 
$*$-homomorphisms. 
If at least one of $A_i$'s is QWEP, then 
the $*$-homomorphism 
$\pi_1\otimes\pi_2\colon C^*\F\otimes C^*\F\to A_1\otimes_{\max}A_2$
is continuous w.r.t.\ the minimal tensor product.
\end{cor}
\begin{proof}
Thanks to Proposition \ref{sep}, 
we may assume that $A_i$'s are separable 
and the free group $\F$ is countable. 
Suppose that $A_2=B/J$ and $B$ has the WEP. 
Since $C^*\F$ has the LLP, $\pi_2$ lifts to a ucp map 
$\psi_2\colon C^*\F\to B$.
It follows that 
$$\pi_1\otimes\pi_2\colon C^*\F\otimes_{\min}C^*\F
\stackrel{\id\otimes\psi_2}{\longrightarrow}
C^*\F\otimes_{\min}B=C^*\F\otimes_{\max}B
\stackrel{\pi_1\otimes Q}{\longrightarrow}
A_1\otimes_{\max}A_2$$
is continuous.
\end{proof}

Choi \cite{choirfd} proved that the $C^*$-algbera $C^*\F_\infty$ 
is residually finite dimensional, i.e., it has a faithful family 
of finite dimensional representations. 
It follows that $C^*\F_\infty$ has a faithful trace. 
Kirchberg \cite{klp} observed that his conjecture is equivalent 
to the above for $\F_\infty\times\F_\infty$. 
We note that Bekka \cite{bekka} proved that 
the full $C^*$-algebra of a residually finite group 
need not be residually finite dimensional (e.g.\ $SL_3(\Z)$). 

\begin{prop}\label{conj2}
The following conjectures are equivalent. 

\begin{enumerate}[\rm (i)]
\item
We have 
$C^*\F_\infty\otimes_{\min}C^*\F_\infty
=C^*\F_\infty\otimes_{\max}C^*\F_\infty$. 

\item
The full $C^*$-algebra $C^*(\F_\infty\times\F_\infty)$ 
is residually finite dimensional. 

\item
The full $C^*$-algebra $C^*(\F_\infty\times\F_\infty)$ 
has a faithful trace. 

\item
The $C^*$-algebra $C^*\F_\infty$ has the WEP. 

\item
All separable $C^*$-algebras are QWEP. 

\item
The LLP implies the WEP. 
\end{enumerate}
\end{prop}
\begin{proof}
We note that 
$C^*(\F_\infty\times\F_\infty)
=C^*\F_\infty\otimes_{\max}C^*\F_\infty$ canonically. 
The implications (i)$\Rightarrow$(ii)$\Rightarrow$(iii) are trivial
and (iii)$\Rightarrow$(i) follows from Lemma \ref{ft} below. 

The implications 
(i)$\Rightarrow$(iv)$\Rightarrow$(v) and 
(vi)$\Rightarrow$(i) follow from Corollary \ref{wepllp}. 
For (v)$\Rightarrow$(vi), we prove that 
a $C^*$-algebra $A$ with the LLP and QWEP has the WEP. 
Let $A\subset\B(\hh)$ and let $\pi$ be a quotient onto $A$ 
from a $C^*$-algebra $B$ with the WEP. 
Fix a finite dimensional operator system $E\subset A$, 
and let $\psi_E\colon E\to B$ be a ucp lifting. 
Since $B$ has the WEP, $\psi_E$ extends to a ucp map 
$\tilde{\psi}_E\colon\B(\hh)\to B^{**}$. 
It follows that $\p_E=\pi^{**}\tilde{\psi}_E$ 
is a ucp map which coincides with the identity on $E$. 
Any cluster point of the net of ucp maps $\p_E\colon\B(\hh)\to A^{**}$ 
in the pointwise weak$^*$ topology is a desired weak expectation. 
This completes the proof. 
\end{proof}

\begin{lem}\label{ft}
Any trace on the maximal tensor product $A_1\otimes_{\max}A_2$ 
factors through the minimal tensor product $A_1\otimes_{\min}A_2$. 
\end{lem}
\begin{proof}
It suffices to show the assertion for extremal traces. 
It is well-known and not too hard to see that 
a trace is extremal if and only if its GNS representation 
generates a finite factor. 
Let $\tau$ be an extremal trace on $A_1\otimes_{\max}A_2$ 
with the GNS representation $\pi$ and 
let $\pi_i\colon A_i\to M:=\pi(A_1\otimes_{\max}A_2)''$ 
be the restriction of $\pi$ to $A_i$. 
It follows that $M_i:=\pi_i(A_i)''$ are commuting 
von Neumann subalgebras of $M$, which have to be factors. 
By uniqueness of the trace on a finite factor, 
we have $\tau_M(a_1a_2)=\tau_{M_1}(a_1)\tau_{M_2}(a_2)$ 
for $a_i\in\pi_i(A_i)$, which means that $M=M_1\bar{\otimes}M_2$.
\end{proof}
Recall the notations used in the proof of Theorem \ref{fb}. 
According to Lemma~\ref{lin}, 
$C^*\F_\infty\otimes_{\max}C^*\F_\infty
=C^*\F_\infty\otimes_{\min}C^*\F_\infty$ 
if and only if the formal identity $\theta_n$ 
from $E_n\otimes_{\min}E_n$ 
into $C^*\F_\infty\otimes_{\max}C^*\F_\infty$ is 
completely contractive for every/some $n\geq 3$. 
It follows from a Grothendieck-type factorization theorem of 
Junge \cite{junge}, Paulsen \cite{pmax} 
and Pisier \cite{pintro} that $\|\theta_n\|\le 2$ for all $n\in\N$. 
However, no nontrivial estimates of $\|\theta_n\|_{\cb}$ is known. 
It is known that the $C^*$-subalgebra in 
$C^*\F_\infty\otimes_{\min}C^*\F_\infty$
generated by its ``diagonal'' $\{ s\otimes s : s\in\F_\infty\}$ 
is canonically $*$-isomorphic to $C^*\F_\infty$. 
Indeed, since $C^*_{\mathrm{red}}\F_\infty$ is QWEP 
(cf.\ the remark at the end of Section 4),
the $*$-homomorphism 
$C^*\F_\infty\otimes_{\min}C^*\F_\infty
\to C^*_{\mathrm{red}}\F_\infty\otimes_{\max}C^*_{\mathrm{red}}\F_\infty$
is continuous thanks to Corollary \ref{lrp}. 
The claim now follows from 
Pisier's observation \cite{pintro} that the ``diagonal'' in 
$C^*_{\mathrm{red}}\G\otimes_{\max}C^*_{\mathrm{red}}\G$ 
is canonically $*$-isomorphic to $C^*\G$ for any discrete group $\G$. 

The following is proved by Pisier \cite{psimple} 
and Boca \cite{boca}. 

\begin{prop}\label{fplp}
The LLP (resp.\ the LP for separable $C^*$-algebras) 
is stable under a full free product. 
\end{prop}

By modifying Pisier's proof \cite{psimple}, 
one can prove that the LLP is stable under 
a full amalgamated free product over 
a finite dimensional $C^*$-subalgebra.
Indeed, if $A_i$ ($i=1,2$) are $C^*$-algebras with 
a common $C^*$-subalgebra $B$, then the 
linear span of $A_1A_2$ in the full amalgamated free product 
$A_1\ast_B A_2$ is canonically completely isometrically isomorphic 
to the relative Haagerup tensor product 
$A_1\otimes_B^hA_2:=A_1\otimes^hA_2/N_B$, 
where $N_B=\mathop{\overline{\mathrm{span}}}
\{a_1x\otimes a_2-a_1\otimes xa_2 : a_i\in A_i,\,x\in B\}$. 
Hence if $B$ is finite dimensional, then there is 
a completely contractive lifting from $A_1\otimes_B^hA_2$ 
into $A_1\otimes^hA_2$. 
It follows from Lemma \ref{lin} that the canonical quotient map 
from $A_1\ast A_2$ onto $A_1\ast_BA_2$ is locally ucp liftable 
whenever $B$ is finite dimensional. 
It is not known whether the full amalgamated free product 
(over a finite dimensional $C^*$-subalgebra) preserves the LP.

The LLP (resp.\ the LP) also passes to a 
(resp.\ separable) weakly cp complemented $C^*$-subalgebra. 
It is not known whether the LLP (resp.\ the LP) is 
stable under the maximal tensor product, or equivalently 
whether the full $C^*$-algebra $C^*(\F_\infty\times\F_\infty)$ 
has the LLP (resp.\ the LP). 
It seems related to the QWEP conjecture, 
but we have found no logical connection. 
Although it seems full group $C^*$-algebras rarely have the LLP, 
there is no example of groups whose full $C^*$-algebra 
is known to fail the LLP. 
It was observed in \cite{univ} that there is a group 
whose full $C^*$-algebra does not have the LP. 
We note that the full $C^*$-algebra $C^*\G$ has 
the LLP if and only if any positive definite function 
from $\G$ into the Calkin algebra $\B(\ell_2)/\K(\ell_2)$ 
has a positive definite lifting. 

\section{Permanence Properties of the QWEP}

Following Kirchberg \cite{klp}, 
we study the permanence properties of the QWEP. 

\begin{prop}\label{perm1}
We have the following. 

\begin{enumerate}[\rm (i)]
\item
If $A_i$ is QWEP for all $i\in I$, then 
so is $\prod_{i\in I}A_i$. 

\item
If $A\subset B$ is weakly cp complemented and $B$ is QWEP, 
then so is $A$. 

\item
Let $(A_i)_{i\in I}$ be an increasing net of (possibly non-unital) 
$C^*$-subalgebras in $A$ (resp. $M$) whose union is dense in norm 
(resp.\ weak$^*$) topology. 
If all $A_i$ are QWEP, then so is $A$ (resp.\ $M$). 

\item
A $C^*$-algebra $A$ is QWEP if and only if 
the second dual $A^{**}$ is QWEP. 

\item
If $A$ is QWEP and $B$ is nuclear, then $A\otimes_{\min}B$ is QWEP. 
If $M$ and $N$ are QWEP, then so is $M\bar{\otimes}N$. 

\item
If $A$ (resp.\ $M$) is QWEP and $\alpha$ is 
an action of amenable group $\G$, then 
$\G\ltimes_\alpha A$ (resp.\ $\G\ltimes_\alpha M$) is QWEP. 

\item
The commutant $M'$ is QWEP if and only if $M$ is QWEP. 

\item
Let $M=\int^\oplus M(\gamma)\,d\gamma$ be 
the direct integral of separable von Neumann algebras. 
Then, $M$ is QWEP if and only if $M(\gamma)$ are QWEP for a.e.\ $\gamma$. 
\end{enumerate}
\end{prop}

\begin{proof}
Ad(i): 
This follows from Lemma \ref{compl}. 

Ad(ii): 
Let $J$ be an ideal in $C$ with the WEP such that $B=C/J$ 
and let $\pi\colon C\to B$ be the quotient map. 
Then, the $C^*$-subalgebra $\pi^{-1}(A)$ is weakly cp complemented 
in $C$ and thus has the WEP. 
Indeed, this follows from Lemma \ref{compl} and the fact 
$$\pi^{-1}(A)^{**}=J^{**}\oplus A^{**}\subset J^{**}\oplus B^{**}=C^{**}.$$ 
Hence, $A$ is a quotient of $\pi^{-1}(A)$ which has the WEP. 

Ad(iii): 
Let $M$ be the weak$^*$-closure of $A=\overline{\bigcup A_i}$. 
We will prove that $M$ is QWEP. 
(Then the QWEP property of $A$ follows from (ii) by considering 
the case where $M=A^{**}$.)
By the Kaplansky density theorem, 
$\ball(A)$ is strong$^*$-dense in $\ball(M)$. 
Hence, amplifying the directed set $I$ if necessary, 
we may assume that 
for any $x\in\ball(M)$, 
there is a net $(a_i)_{i\in I}\in\ball(\prod_{i\in I}A_i)$ 
such that $x=\mbox{strong$^*$-}\lim_{i\in I}a_i$
We denote by $C=\prod_{i\in I}A_i$. 
The $C^*$-algebra $C$ is QWEP by (i). 
Let 
$$B=\{ (a_i)_{i\in I}\in C : 
\mbox{strong$^*$-}\lim_{i\in I}a_i\mbox{ exists in }M\}$$
and let $\pi\colon B\to M$ be the map which takes 
$(a_i)_{i\in I}$ to its limit. 
Since the adjoint-operation and product is jointly 
strong$^*$-continuous on bounded sets, 
$B$ is a $C^*$-subalgebra of $C$ and 
$\pi$ is a surjective $*$-homomorphism onto $M$. 
We claim that $B$ is weakly cp complemented in $C$. 
Let $(e_j)_j$ be an increasing approximate unit for $J=\ker\pi$ 
and let $e=\lim_j e_j\in C^{**}$. 
We note that $e$ is in the center of $B^{**}$ and 
$B^{**}=B^{**}e\oplus B^{**}(1-e)=J^{**}\oplus M^{**}$. 
Since $J$ is hereditary in $C$ 
(i.e., $e_ixe_i\in J$ for all $x\in C$), 
we have $eC^{**}e=J^{**}$. 
For a fixed free ultrafilter $\omega$ on $I$, the map 
defined by 
$$\tilde{\pi}\colon C\ni(a_i)_{i\in I}
\longmapsto\mbox{weak$^*$-}\lim_{i\in\omega}a_i\in M$$ 
is a ucp extension of $\pi$. 
This ucp extension gives rise to a conditional expectation 
from $(1-e)C^{**}(1-e)$ onto $B^{**}(1-e)=M^{**}$. 
It follows that 
$$B^{**}=J^{**}\oplus M^{**}\subset eC^{**}e\oplus(1-e)C^{**}(1-e)$$
is cp complemented. 
This proves our claim and we are done by (ii).

Ad(iv): 
This follows from (ii) and (iii). 

\noindent\textit{Note.} 
We are not going to use (v) to (vii) of this proposition 
and the proof requires some results which 
will be proved later. 

Ad(v): 
The first assertion follows from the fact that 
the tensor product of a WEP $C^*$-algebra with a nuclear $C^*$-algebra 
again has the WEP. 
The second assertion follows from Corollary \ref{frep} 
(or Corollary \ref{emb}) and the following result of Tomiyama \cite{tomiyama}. 
If $\p\colon M_1\to M_2$ is a (not necessarily weak$^*$-continuous) 
ucp map between von Neumann algebras, then 
the map $\p\otimes\id_N\colon M_1\otimes N\to M_2\otimes N$ 
extends to a ucp map 
$\p\bar{\otimes}\id_N\colon M_1\bar{\otimes}N\to M_2\bar{\otimes}N$. 
Indeed, this follows from the fact that $(M_2)_*\otimes N_*$ 
is dense in $M_2\bar{\otimes}N$.

Ad(vi): 
The assertion for von Neumann algebras follows 
from that for $C^*$-algebras. 
Let $(A,\G,\alpha)$ be a $C^*$-dynamical system with $\G$ amenable. 
We claim that there are nets of ucp maps 
$\psi_i\colon\G\ltimes A\to\M_{n(i)}(A)$ 
and $\p_i\colon \M_{n(i)}(A)\to\G\ltimes A$ 
such that $\p_i\psi_i$ converges to $\id_{\G\ltimes A}$ pointwisely. 
Once this is proved, then the QWEP property of $\G\ltimes A$ 
follows from that of $A$ (cf.\ Corollary \ref{frep}). 

Take a faithful representation $\pi\colon A\to\B(\hh)$ 
with a unitary action $u$ of $\G$ on $\hh$ 
which implements the action $\alpha$. i.e., 
$\pi(\alpha_s(a))=\Ad u_s(\pi(a))$ for $a\in A$. 
Then, we have 
$\G\ltimes A\subset\B(\ell_2\G\otimes\hh)$ 
where $\G\ni s\mapsto\lambda_s\otimes 1\in\B(\ell_2\G\otimes\hh)$ and 
$A\ni a\mapsto\tilde{\pi}(a)\in\B(\ell_2\G\otimes\hh)$ with 
$\tilde{\pi}(a)(\delta_s\otimes\xi)
=\delta_s\otimes\alpha_{s^{-1}}(a)\xi$. 
We fix a finite subset $F\subset\G$. 
The compression of $\B(\ell_2\G)$ to $\B(\ell_2F)=\M_F$ defines 
a ucp map $\psi_F\colon\G\ltimes A\to\M_F(A)$. 
Let $\p_F\colon\M_F(A)\to\B(\ell_2\G\otimes\hh)$ be the ucp map 
defined by $\p_F(x)=V_F^*(1\otimes x)V_F$, where 
the isometry $V_F$ is given by 
$$V_F\colon\ell_2\G\otimes\hh\ni\delta_s\otimes\xi
\longmapsto \frac{1}{\sqrt{|F|}}\sum_{t\in F}
\delta_{t^{-1}s}\otimes\delta_t\otimes u_{t^{-1}s}\xi
\in\ell_2\G\otimes\ell_2F\otimes\hh.$$
A direct computation shows that 
$\p_F(e_{st,t}\otimes a)
=\frac{1}{|F|}(\lambda_s\otimes1)\tilde{\pi}(\alpha_t(a))$. 
Therefore, $\p_F$ maps into $\G\ltimes A$ and 
$\p_F\psi_F(\lambda_s\tilde{\pi}(a))
=\frac{|sF\cap F|}{|F|}\lambda_s\tilde{\pi}(a)$.
This proves the claim.

Ad(vii): 
This follows from standard representation theory 
of von Neumann algebras.

Ad(viii): 
Suppose $M$ is QWEP and let $A$ be a separable weak$^*$-dense 
$C^*$-subalgebra which is weakly cp complemented in $M$ 
(cf.\ Proposition \ref{sep}).
Then, $A$ is QWEP by (ii) in this proposition. 
Since almost all $M(\gamma)$ arise as the weak$^*$-closure 
of representations of $A$, they are QWEP.
Conversely, suppose $M(\gamma)$ are QWEP for a.e.\ $\gamma$ 
and let $\pi$ (resp.\ $\pi'$) be a $*$-homomorphism from $C^*\F_\infty$ 
onto a weak$^*$-dense $C^*$-subalgebra of $M$ (resp.\ $M'$). 
Since $\pi_\gamma\cdot\pi'_\gamma$ is coninuous w.r.t.\ 
the minimal tensor product for a.e.\ $\gamma$ (cf.\ Corollary \ref{lrp}), 
the representation 
$\pi\cdot\pi'=\int^\oplus\pi_\gamma\cdot\pi'_\gamma\,d\gamma$ 
is also coninuous w.r.t.\ the minimal tensor product.
If $C^*\F_\infty\subset\B(\hh)$, 
then one can show that 
$\pi\colon C^*\F_\infty\to M$ extends to a ucp map 
$\tilde{\pi}\colon\B(\hh)\to M$ 
(cf.\ the proof of Theorem \ref{lance}). 
The QWEP property of $M$ now follows from Corollary \ref{frep}. 
\end{proof}

It is not known whether the minimal tensor product of 
QWEP $C^*$-algebras is again QWEP, or equivalently 
whether $\B(\ell_2)\otimes_{\min}\B(\ell_2)$ is QWEP or not.  
We note that this $C^*$-algebra fails the WEP as shown in \cite{wep}. 

Haagerup \cite{hmap} showed that there is a sequence 
of (complete) contractions 
$\p_n\colon L\F_2\to C^*_{\mathrm{red}}\F_2=:A$ 
such that $\lim_n\p_n(a)=a$ for $a\in A$. 
Let $\p\colon L\F_2\to A^{**}$ 
be a cluster point of the sequence $\p_n$. 
Then, $\p$ is a unital (complete) contraction such that $\p|_A=\id_A$. 
This shows that $C^*_{\mathrm{red}}\F_2$ is weakly cp complemented 
in $L\F_2$ (cf.\ Theorem \ref{tomiyama}). 
(It follows that $C^*_{\mathrm{red}}\F_2$ 
is weakly cp complemented in any $C^*$-superalgebra 
with a trace.) 
Since $L\F_2$ is QWEP (cf.\ Proposition \ref{rf}), 
the reduced group $C^*$-algebra $C^*_{\mathrm{red}}\F_2$ 
is also QWEP. 
However, it seems unknown whether the reduced group $C^*$-algebra 
$C^*_{\mathrm{red}}\G$ of, say, $\G=SL_3(\Z)$ is QWEP or not. 

\section{Finite Representability in the Trace Class}

Following Kirchberg \cite{klp}, 
we study the relation between the QWEP property 
and finite representability.

\begin{defn}
The Banach-Mazur distance $\dist$ between two 
Banach spaces $E$ and $F$ is defined by 
$$\dist(E,F)=\inf\{ \|\theta\|\,\|\theta^{-1}\| : 
\theta\colon E\to F\mbox{ a continuous linear isomorphism}\}.$$
We put $\dist(E,F)=\infty$ when $E$ is not isomorphic to $F$. 
We say a Banach space $X$ is finitely representable in 
a Banach space $Y$ if for every finite dimensional subspace $E$ of $X$, 
one has $\inf\{\dist(E,F): F\subset Y\}=1$. 

Similarly, the cb Banach-Mazur distance $\dist_{\cb}$ 
and os finite representability are defined by just replacing 
the usual norm with the cb norm in the above definitions. 
\end{defn}

For a Banach space $X$, we denote by $\ball(X)$ (resp. $\cball(X)$) 
the open (resp.\ closed) unit ball of $X$. 
A map $\p\colon X\to Y$ is called a \emph{metric surjection}
if it maps $\ball(X)$ onto $\ball(Y)$, 
i.e., $X/\ker\p=Y$ isometrically. 
If $Y$ embeds into $X$ isometrically, 
then the transpose of the isometric inclusion is 
a weak$^*$-continuous metric surjection which maps 
$\cball(X^*)$ onto $\cball(Y^*)$. 

\begin{prop}\label{ms}
For von Neumann algebras $M$ and $N$, we have the following. 

\begin{enumerate}[\rm (i)]
\item
If $M_*$ embeds into $N_*$ isometrically, 
then there is a weak$^*$-continuous 
unital metric surjection $\p$ from $N$ onto $M$. 

\item
If there is a weak$^*$-continuous unital metric surjection $\p$ 
from $N$ onto $M$, then there are projections $p$ and $q$ in $N$ 
such that $M$ is $*$-isomorphic to a (normally)
cp complemented von Neumann subalgebra of $pNp\oplus(qNq)^{\op}$.
\end{enumerate}
\end{prop}
\begin{proof}
Ad(i): 
Let $\p_0\colon N\to M$ be the transpose of 
an isometric embedding of $M_*$ into $N_*$. 
Let $v$ be an extreme point of the weak$^*$-closed face  
$\{ x\in N : \| x\|\le 1\mbox{ and }\p_0(x)=1\}$ 
of $\cball(N)$. 
It follows from a standard argument (cf.\ \cite{sakai}) that 
$v$ is a partial isometry and $\p_0(vv^*x)=\p_0(x)$. 
Hence, the map $\p\colon N\to M$ defined by $\p(x)=\p_0(vx)$ 
is the desired unital metric surjection. 

Ad(ii): 
Let $\p$ be a weak$^*$-continuous 
unital metric surjection from $N$ onto $M$.
Let $e\in N$ be the support projection of $\p$, i.e., 
$1-e=\sup\{ a\in N : 0\le a\le 1\mbox{ and }\p(a)=0\}$. 
(We note that $a\geq0$ and $\p(a)=0$ implies that 
$\p$ is zero on the support projection of $a$.)
It follows that $\p(a)=\p(eae)$ for every $a\in N$ and 
that $\p(a)\neq 0$ for all positive non-zero element $a\in eNe$.
Let $C\subset eNe$ be the Jordan multiplicative domain for $\p|_{eNe}$, 
i.e., 
\begin{align*}
C &=\{ a\in eNe : \p(a^*\circ a)=\p(a)^*\circ\p(a)\}\\
 &=\{ a\in eNe : \p(a\circ x)=\p(a)\circ\p(x)\mbox{ for all }x\in eNe\}.
\end{align*}
By Corollary \ref{mdom} and the following remark, 
$C$ is a weak$^*$-closed Jordan subalgebra in $eNe$ and 
$\p|_C$ is a normal Jordan isomorphism from $C$ onto $M$. 
Hence, $\theta=(\p|_C)^{-1}$ is a normal Jordan isomorphism 
from $M$ onto $C$ with $\p\theta=\id_M$. 
By Theorem \ref{jordan}, 
there is a projection $p$ in $C'\cap eNe$ with $q:=e-p$ such that 
the map $\theta_1\colon M\ni x\mapsto\theta(x)p\in pNp$ 
(resp.\ $\theta_2\colon M\ni x\mapsto\theta(x)q\in qNq$) 
is a $*$-homomorphism (resp. $*$-antihomomorphism). 
It follows that the map 
$$\tilde{\theta}\colon M\ni x
\longmapsto(\theta_1(x),\theta_2(x)^{\op}\big)
\in pNp\oplus(qNq)^{\op}$$
defines a normal unital injective $*$-homomorphism 
and the map 
$$\tilde{\p}\colon pNp\oplus(qNq)^{\op}
\ni (a,b^{\op})\longmapsto\p(a+b)\in M$$
defines a unital positive map such that $\tilde{\p}\tilde{\theta}=\id_M$. 
Therefore, the map $\tilde{\theta}\tilde{\p}$ 
is a weak$^*$-continuous contractive projection onto 
the von Neumann subalgebra 
$\tilde{\theta}(M)$ in $pNp\oplus(qNq)^{\op}$. 
We complete the proof by Theorem \ref{tomiyama}. 
\end{proof}

\begin{cor}\label{frep}
For a von Neumann algebra $M$, the following are equivalent. 

\begin{enumerate}[\rm (i)]
\item
The von Neumann algebra $M$ is QWEP. 

\item
There are a Hilbert space $\hh$ and a ucp map 
$\p\colon\B(\hh)\to M$ which maps $\cball(\B(\hh))$ 
onto $\cball(M)$. 

\item
The predual $M_*$ of $M$ is (os) finitely representable 
in the trace class $S_1$. 

\item 
There is a QWEP $C^*$-algebra $A$ and 
a contraction $\p\colon A\to M$ 
such that $\p(\cball(A))$ is weak$^*$-dense in $\cball(M)$. 
\end{enumerate}
Moreover, we can choose the Hilbert space $\hh$ in $\mathrm{(ii)}$ 
separable when so is $M$. 
\end{cor}

\begin{proof}
Ad(i)$\Rightarrow$(ii): 
Let $B\subset\B(\hh)$ be a $C^*$-algebra with the WEP 
and let $\pi$ be a surjective $*$-homomorphism from $B$ onto $M$. 
Since $B$ has the WEP, there is a ucp map 
$\psi\colon\B(\hh)\to B^{**}$ with $\psi|_B=\id_B$. 
It follows that the composition $\p=\tilde{\pi}\psi$ of 
$\psi$ and the normal extension $\tilde{\pi}\colon B^{**}\to M$ 
of $\pi$ is a ucp map such that $\p(\cball(\B(\hh)))=\cball(M)$. 
If $M$ is separable, then there is a separable $C^*$-subalgebra $A$ 
in $\B(\hh)$ such that $\p(\cball(A))$ is weak$^*$-dense in $\cball(M)$. 
Let $\theta\colon A\to\B(\ell_2)$ be a faithful representation 
and let $\sigma\colon\B(\ell_2)\to\B(\hh)$ be 
a ucp extension of $\theta^{-1}$. 
It follows that $\p_0=\p\sigma$ is a ucp map such that 
$\p_0(\cball(\B(\ell_2)))$ is weak$^*$-dense in $\cball(M)$. 
Fixing a free ultrafilter $\omega$ on $\N$, 
we define $\Phi\colon\prod_{n\in\N}\B(\ell_2)\to M$ by 
$$\Phi((a_n)_{n\in\N})
=\mbox{weak$^*$-}\lim_{n\in\omega}\p_0(a_n)\in M.$$
Since $M_*$ is separable, the ucp map $\Phi$ maps 
$\cball(\prod_{n\in\N}\B(\ell_2))$ onto $\cball(M)$. 

Ad(ii)$\Rightarrow$(iii): 
A metric surjection between $C^*$-algebras which is ucp 
is automatically a complete metric surjection 
(cf.\ Corollary \ref{mdom}). 
If $\p\colon\B(\hh)\to M$ is a ucp metric surjection, 
then $\p^*\colon M^*\to\B(\hh)^*$ is completely isometric. 
By the principle of local reflexivity, 
$\B(\hh)^*$ is finitely representable in $\B(\hh)_*=S_1(\hh)$. 
The operator space analogue of the principle of local reflexivity 
for $S_1(\hh)$ is a deep theorem of Junge \cite{junge} 
(see \cite{ejr} for a more general result). 

Ad(iii)$\Rightarrow$(iv): 
Let $(E_i)_{i\in I}$ be an increasing net of 
finite dimensional subspaces in $M_*$ with $\bigcup E_i=M_*$. 
By the assumption, 
there is an embedding $\psi_i\colon E_i\to S_1$ for each $i\in I$ 
such that $\|\psi_i\|\le1+(\dim E_i)^{-1}$ and $\|\psi_i^{-1}\|\le1$.
Fixing a free ultrafilter $\omega$ on $I$, we define a 
contraction $\p\colon\prod_{i\in I}\B(\ell_2)\to M$ by 
$$\ip{\p((a_i)_{i\in I}),f}=\lim_{i\in\omega}\ip{x_i,\psi_i(f)}$$
for $(a_i)_{i\in I}\in\prod_{i\in I}\B(\ell_2)$ and $f\in M_*$ 
(since $f\in E_i$ eventually, 
the limit in the above definition makes sense).  
We give ourselves an arbitrary $x\in\cball(M)$. 
For each $i\in I$, let $a_i\in\B(\ell_2)=S_1^*$ 
be a Hahn-Banach extension of $x\psi_i^{-1}\in\psi_i(E_i)^*$. 
It follows that $a=(a_i)_{i\in I}\in\cball(\prod_{i\in I}\B(\ell_2))$ 
is such that $\p(a)=x$. 

Ad(iv)$\Rightarrow$(i): 
By the Kaplansky density theorem, 
there is an directed set $I$ such that 
for any $x\in\ball(M)$, 
there is a net $(a_i)_{i\in I}$ in $\ball(A)$ 
such that $x=\mbox{weak$^*$-}\lim_{i\in I}\p(a_i)$
It follows that the map 
$$\prod_{i\in I}A\ni(a_i)_{i\in I}\longmapsto
\mbox{weak$^*$-}\lim_{i\in\omega}\p(a_i)\in M$$
defines a metric surjection from $\prod_{i\in I}A$ onto $M$. 
Taking the transpose, we obtain an isometric embedding 
of $M_*$ into $(\prod_{i\in I}A)^*$. 
By Proposition \ref{ms}, $M$ is $*$-isomorphic to 
a (not necessarily unital) cp complemented von Neumann subalgebra 
of $N\oplus N^{\op}$, where $N=(\prod_{i\in I}A)^{**}$. 
We conclude by Proposition \ref{perm1} that $M$ is QWEP. 
\end{proof}

%
%
%
\section{Connes' Embedding Problem}

Let $\omega$ be a fixed free ultrafilter on $\N$ and 
let $R$ be the hyperfinite $\mathrm{II}_1$-factor whose
faithful normal trace is denoted simply by $\tau$. 
We define a trace $\tau_\omega$ on $\prod_{n\in\N}R$ by 
$\tau_\omega((a_n))=\lim_{n\to\omega}\tau(a_n)$. 
The ultrapower $R^\omega$ of the hyperfinite $\mathrm{II}_1$-factor $R$ 
is defined to be the von Neumann algebra generated by 
the GNS representation associated with $\tau_\omega$. 
We still denote by $\tau_\omega$ the faithful trace on $R^\omega$. 
It turns out that 
$N_\omega=\{a\in\prod R : \tau_\omega(a^*a)=0\}$ 
is a maximal ideal in $\prod R$ and 
$R^\omega=(\prod R)/N_\omega$ (as a $C^*$-algebra) 
is a $\mathrm{II}_1$-factor which is not separable. 

Let $A$ be a $C^*$-algebra and $\tau$ be a trace on $A$ 
with the GNS-triplet $(\pi_\tau,\hh_\tau,\xi_\tau)$.
We denote by $\pi_\tau^c$ the representation 
of the conjugate $C^*$-algebra $\bar{A}$ of $A$ on $\hh_\tau$
defined by 
$\pi_\tau^c(\bar{b})\pi_\tau(a)\xi_\tau=\pi_\tau(ab^*)\xi_\tau$ 
for $a, b\in A$. 
(The conjugate $C^*$-algebra $\bar{A}$ 
is $*$-isomorphic to the opposite $C^*$-algebra $A^{\op}$ 
via $\bar{A}\ni\bar{a}\mapsto (a^*)^{\op}\in A^{\op}$.)
This gives rise to a representation $\sigma_\tau$
of $A\otimes\bar{A}$ on $\hh_\tau$ given by 
$$\sigma_\tau\colon A\otimes\bar{A}
 \ni\sum_ka_k\otimes\bar{b_k}
 \longmapsto\sum_k\pi_\tau(a_k)\pi_\tau^c(\bar{b}_k)\in\B(\hh_\tau)$$ 
which is continuous w.r.t.\ the maximal tensor norm. 
It follows that the linear functional 
$\mu_\tau$ on $A\otimes\bar{A}$ 
given by 
$$\mu_\tau\colon A\otimes\bar{A}
 \ni\sum_ka_k\otimes\bar{b}_k
 \longmapsto\tau(\sum_ka_kb_k^*)\in\C$$ 
is also continuous w.r.t.\ the maximal tensor norm. 
Connes showed that a $\mathrm{II}_1$-factor $(A,\tau)$
is injective iff $\sigma_\tau$ (or $\mu_\tau$) is continuous 
w.r.t.\ the minimal tensor norm on $A\otimes\bar{A}$. 
The following theorem of Kirchberg generalizes Connes' characterization. 

\begin{thm}\label{trace}
For a trace $\tau$ on a $C^*$-algebra $A$ in $\B(\hh)$, 
the following are equivalent. 

\begin{enumerate}[\rm (i)]
\item
The trace $\tau$ extends to an $A$-central state $\p$ on $\B(\hh)$,
i.e., the trace $\tau$ extends to a state $\p$ on $\B(\hh)$ such that 
$\p(ax)=\p(xa)$ for every $a\in A$ and $x\in\B(\hh)$. 

\item
There is a net of ucp maps $\theta_i\colon A\to\M_{n(i)}$ 
such that $\tau(a)=\lim_i\tr_{n(i)}(\theta_i(a))$ 
and $\lim_i\tr_{n(i)}(\theta_i(ab^*)-\theta_i(a)\theta_i(b)^*)=0$ 
for every $a,b$ in $A$. 

\item[\rm (ii')]
The trace $\tau$ is liftable, i.e., 
there is a $*$-homomorphism $\theta\colon A\to R^\omega$ 
with a ucp lifting $\widehat{\theta}\colon A\to\prod R$ 
such that $\tau=\tau_\omega\theta$. 

\item
The functional $\mu_\tau$ is continuous 
w.r.t.\ the minimal tensor norm on $A\otimes\bar{A}$. 

\item
The representation $\sigma_\tau$ is continuous 
w.r.t.\ the minimal tensor norm on $A\otimes\bar{A}$. 
\end{enumerate}
Moreover if $A\cap\K(\hh)=\{0\}$, then one can choose 
the ucp maps $\theta_i\colon A\to\M_{n(i)}$ 
in the condition $\mathrm{(ii)}$ as compressions 
to finite dimensional subspaces of $\hh$. 
\end{thm}
\begin{proof}
We will prove (i)$\Rightarrow$(ii)$\Rightarrow$(iii)$\Rightarrow$(iv)
$\Rightarrow$(i) and (ii)$\Rightarrow$(ii')$\Rightarrow$(i).

Ad (i)$\Rightarrow$(ii): 
The proof is taken from \cite{haagerup}. 
To prove (ii), we give ourselves a finite set $\FF$ 
of unitaries in $A$ and $\e>0$. 
We approximate the $A$-central state $\p$ by $\Tr(h\,\cdot\,)$, 
where $h$ is a positive trace class operator with $\Tr(h)=1$. 
By a standard approximation argument, we find such $h$ that 
$|\Tr(hu)-\tau(u)|<\e$ and $\|h-uhu^*\|_{1,\Tr}<\e$ 
for $u\in\FF$. 
Further, we may assume that $h$ is of finite rank and 
has no irrational eigenvalues; 
let $p_1/q,\ldots,p_m/q$ ($p_1,\ldots,p_m\in\N$ and $q=\sum_kp_k$) 
be the non-zero eigenvalues of $h$ with the corresponding 
eigenvectors $\zeta_1,\ldots,\zeta_m\in\hh$. 
Put $p=\max\{p_1,\ldots,p_m\}$ and define an isometry
$V_k\colon\ell_2^{p_k}\to\hh\otimes\ell_2^p$ for each $k$ 
by $V_k\delta_i=\zeta_k\otimes\delta_i$ for $i=1,\ldots,p_k$. 
Since $V_k$'s have orthogonal ranges, the sum of $V_k$'s 
gives rise to an isometry 
$V\colon\bigoplus_{k=1}^m\ell_2^{p_k}\to\hh\otimes\ell_2^p$. 
Identifying $\M_q$ with $\B(\bigoplus_{k=1}^m\ell_2^{p_k})$, 
we obtain a ucp map $\theta\colon A\to\M_q$ defined by 
$\theta(a)=V^*(a\otimes 1)V$. 
It follows that $\tr_q\theta(a)=\Tr(ha)$ for $a\in A$ and 
that, by denoting $u_{k,l}=(u\zeta_l,\zeta_k)$, we have 
\begin{align*}
|\Tr(h^{1/2}uh^{1/2}u^*) &- \tr_q(\theta(u)\theta(u)^*)| 
 = \sum_{k,l}
  |u_{k,l}|^2\left((p_kp_l)^{1/2}-\min\{p_k,p_l\}\right)/q\\
 &\le \sum_{k,l}|u_{k,l}|^2 p_k^{1/2}|p_k^{1/2}-p_l^{1/2}|/q\\
 &\le \left(\sum_{k,l}|u_{k,l}|^2p_k/q\right)^{1/2}
  \left(\sum_{k,l}|u_{k,l}|^2(p_k^{1/2}-p_l^{1/2})^2/q\right)^{1/2}\\
 &= \| h^{1/2}u\|_{2,\Tr}\,\|h^{1/2}u-uh^{1/2}\|_{2,\Tr}.
\end{align*}
Recall the Powers-St{\o}rmer inequality \cite{ps} that 
$\|h^{1/2}u-uh^{1/2}\|_{2,\Tr}\le\|h-uhu^*\|_{1,\Tr}^{1/2}$. 
It follows that 
\begin{align*}
\tr_q & (\theta(uu^*) - \theta(u)\theta(u)^*)\\
&= |\Tr(uhu^*)-\Tr(h^{1/2}uh^{1/2}u^*)|
 + |\Tr(h^{1/2}uh^{1/2}u^*) - \tr_q(\theta(u)\theta(u)^*)|\\
&\le 2\|h^{1/2}u-uh^{1/2}\|_{2,\Tr} \le 2\e^{1/2}
\end{align*}
for $u\in \FF$. 
Finally, we note that 
\[
|\tr_q(\theta(ab^*)-\theta(a)\theta(b)^*)|
\le \tr_q(\theta(aa^*)-\theta(a)\theta(a)^*)^{1/2}
  \tr_q(\theta(bb^*)-\theta(b)\theta(b)^*)^{1/2}
\]
for any $a, b\in A$. 
If $A\cap\K(\hh)=\{0\}$, then by Glimm's theorem 
the ucp map $\theta\colon A\to\M_q$ is approximated by 
compressions to $q$-dimensional subspaces in $\hh$. 
Since $A$ is spanned by unitaries, this completes the proof 
of (i)$\Rightarrow$(ii). 

Ad(ii)$\Rightarrow$(iii):
Since 
$\mu_n\colon\M_n\otimes_{\min}\bar{\M}_n
\ni\sum_kx_k\otimes\bar{y}_k\mapsto\tr_n(\sum_kx_ky_k^*)\in\C$ 
are states for all $n$, the net of states 
$\mu_{n(i)}(\theta_i\otimes\bar{\theta_i})$ 
on $A\otimes_{\min}\bar{A}$ is well-defined and 
converges to the functional $\mu_\tau$. 

Ad(iii)$\Rightarrow$(iv):
This follows from the fact that $\xi_\tau$ 
is cyclic for $\sigma_\tau(A\otimes\bar{A})$ 
and the corresponding vector state $\mu_\tau$ 
is continuous w.r.t.\ the minimal tensor product. 

Ad(iv)$\Rightarrow$(i):
Let $\Psi\colon \B(\hh)\otimes_{\min}\bar{A}\to\B(\hh_\tau)$ 
be a ucp extension of $\sigma_\tau$ and let  
$\psi\colon\B(\hh)\ni x\mapsto\Psi(x\otimes 1)\in\B(\hh_\tau)$.
Since $\Psi$ is an $(A\otimes_{\min}\bar{A})$-bimodule map, 
we have $\pi_\tau^c(\bar{b})\psi(x)=\Psi(x\otimes\bar{b})
=\psi(x)\pi_\tau^c(\bar{b})$ 
for all $x\in\B(\hh)$ and $\bar{b}\in\bar{A}$. 
It follows that 
$\psi(\B(\hh))\subset\pi_\tau^c(\bar{A})'=\pi_\tau(A)''$ 
and $\psi|_A=\pi_\tau$ 
(in particular, $\psi$ is an $A$-bimodule map). 
Thus, the state $\p$, given by $\p(x)=(\psi(x)\xi_\tau,\xi_\tau)$ 
for $x\in\B(\hh)$, is a desired $A$-central extension of $\tau$. 

Ad(ii)$\Rightarrow$(ii')$\Rightarrow$(i):
The implication (ii)$\Rightarrow$(ii') follows from 
the definition of $R^\omega$ and the fact that 
$\M_n$'s are isomorphic to subfactors of $R$ for all $n$. 
Now assume (ii') and let $\theta$ and $\widehat{\theta}$ be 
as in the condition. 
Since $\prod R$ is injective, $\widehat{\theta}$ extends to 
a ucp map from $\B(\hh)$ into $\prod R$. 
Composing this with the quotient map $\prod R\to R^\omega$, 
we obtain a ucp extension $\psi\colon\B(\hh)\to R^\omega$ of $\theta$. 
Since $\psi$ is an $A$-bimodule map, 
the state $\p$ given by $\p=\tau_\omega\psi$ is 
a desired $A$-central extension of $\tau$. 
\end{proof}

\begin{cor}\label{emb}
Let $M$ be a separable finite von Neumann algebra. 
Then $M$ is QWEP if and only if $M$ is $*$-isomorphic 
to a von Neumann subalgebra of $R^\omega$. 
\end{cor}
\begin{proof}
Any von Neumann subalgebra of $R^\omega$ is the range of 
a conditional expectation and a fortiori is QWEP. 
Now, assume that $M$ is QWEP and 
let $\tau$ be a normal faithful state on $M$. 
Let $\pi\colon C^*\F_\infty\to M$ be a $*$-homomorphism 
with a weak$^*$-dense range. 
Since $M$ is QWEP, by Corollary \ref{lrp}, the $*$-homomorphism
$\pi\otimes\bar{\pi}\colon 
C^*\F_\infty\otimes_{\min}\overline{C^*\F_\infty}
\to M\otimes_{\max}\bar{M}$ is continuous. 
It follows that Theorem \ref{trace} is applicable to 
the trace $\tau\pi$ on $C^*\F_\infty$; 
there is a $*$-homomorphism 
$\theta\colon C^*\F_\infty\to R^\omega$ 
such that $\tau\pi=\tau_\omega\theta$. 
This means that $M$ is $*$-isomorphic to the von Neumann 
subalgebra generated by $\theta(C^*\F_\infty)$ in $R^\omega$. 
\end{proof}

\section{Groups with the Factorization Property}

Otherwise stated, all groups denoted by the symbol $\G$ are 
assumed countable and discrete.

If Connes' embedding problem has a negative answer, 
then it is quite natural to seek a counterexample 
in the group von Neumann algebras. 
It turns out \cite{kf}, \cite{radulescu2} that 
embeddability of $L\G$ into $R^\omega$ is equivalent to 
that of $\G$ into the unitary group $\U(R^\omega)$ of $R^\omega$.
We say a group $\G$ is hyperlinear if it embeds into $\U(R^\omega)$. 

\begin{prop}
A group $\G$ is hyperlinear if and only if 
$L\G$ is $*$-isomorphic to a von Neumann subalgebra of $R^\omega$. 
\end{prop}
\begin{proof}
One direction is obvious. 
Let $\theta\colon\G\to\U(R^\omega)$ be an injective homomorphism 
and let $\Theta=(\theta_n)_{n\in\N}\colon\G\to\U(\prod_{n\in\N}R)$ 
be any lift, i.e., $\pi\Theta(s)=\theta(s)$ for $s\in\G$. 
We note that $\theta_n$ is asymptotically multiplicative, 
but no longer multiplicative. 
For each $n$, we define 
$$\tilde{\theta}_n\colon\G\ni s
\longmapsto\left[\begin{array}{cc}
\theta_n(s)\otimes 1 & \\ & \theta_n(s)\otimes\theta_n(s)
\end{array}\right]\in\U(\M_2(R\bar{\otimes}R)).$$
Since $|z+z^2|/2<1$ for $z\in\C$ with $|z|\le1$ and $z\neq1$, 
we have $\lim_{n\to\omega}|\tau\tilde{\theta}_n(s)|<1$ 
for all $s\in\G$ with $s\neq1$. 
It follows that for an appropriately chosen 
sequence $k(n)$ of integers, the sequence of functions
$$\psi_n:=\tilde{\theta}_n^{\otimes k(n)}
\colon\G\to\U(\M_2(R)^{\otimes k(n)})\cong\U(R)$$ 
is asymptotically multiplicative and 
$\lim_{n\in\omega}\tau\psi_n(s)=0$ for all $s\in\G$ with $s\neq1$.
This means that for the homomorphism $\psi=\lim\psi_n$, 
the von Neumann subalgebra of $R^\omega$ generated by $\psi(\G)$ 
is canonically $*$-isomorphic to $L\G$. 
\end{proof}

\begin{defn}
We say a group $\G$ has the property $\mathrm{(F)}$ 
(the factorization property) 
if the trace $\tau$ on the full $C^*$-algebra $C^*\G$, 
defined by $\tau(s)=\delta_{1,s}$ for $s\in\G$, is liftable.
\end{defn}

By Theorem \ref{trace}, a group $\G$ has the property $\mathrm{(F)}$
if and only if the representation 
$$\sigma_\G\colon C^*\G\otimes C^*\G\ni \sum_ka_k\otimes b_k
\longmapsto \sum_k\lambda_\G(a_k)\rho_\G(b_k)\in\B(\ell_2\G)$$ 
is continuous w.r.t.\ the minimal tensor norm, where 
$\lambda_\G$ (reps.\ $\rho_\G$) is the left (resp.\ right) 
regular representation of $\G$ on $\ell_2\G$. 
A group $\G$ with the property $\mathrm{(F)}$ is hyperlinear, 
and to the best of our knowledge, all groups known to be 
hyperlinear satisfy the property $\mathrm{(F)}$. 
Let $\G$ be a hyperlinear group and let $N$ be a normal subgroup 
of $\F_\infty$ such that $\G=\F_\infty/N$. 
Then, because of the LLP of $C^*\F_\infty$, 
the trace $\tau_\G$ on $C^*\F_\infty$, 
given by the characteristic function of $N$, 
is liftable. 
Hence, the study on permanence properties of liftable traces 
is more general than that of hyperlinearity.  

If P is a property of groups (e.g.\ being finite) 
then a group $\G$ is said to be residually P if, 
for any $s\in\G$ with $s\neq1$, there is a normal subgroup 
$\Delta$ in $\G$ such that $s\notin\Delta$ and $\G/\Delta$ satisfy P. 
When $\Delta$ is a subgroup of $\G$, the quotient homogeneous space 
$\G/\Delta$ is said to be amenable if there is a $\G$-invariant mean 
on $\ell_\infty(\G/\Delta)$. 
This is equivalent to the existence of a F{\o}lner sequence 
in $\G/\Delta$ (see \cite{greenleaf}). 
As it was observed by Wassermann \cite{wassermann} and 
Kirchberg \cite{kf}, residually finite groups 
and in particular the free groups have the property $\mathrm{(F)}$. 
This also follows from a more general result of Brown 
and Dykema \cite{bd} that the property $\mathrm{(F)}$ is preserved 
under a free product (with amalgamation over a finite subgroup). 

\begin{prop}\label{rf}
A group $\G$ which is residually $\mathrm{(F)}$ 
has the property $\mathrm{(F)}$. 
In particular, residually finite groups have the property $\mathrm{(F)}$. 

Moreover, the property $\mathrm{(F)}$ is preserved under 
taking a subgroup, a supergroup with an amenable quotient, 
a direct product, an increasing union, 
and an amalgamated free product over a finite group. 
\end{prop}
\begin{proof}
The first assertion follows from the fact that 
the set of liftable traces is closed in the pointwise topology. 
We only prove nontrivial claims in the second. 

The claim on a supergroup follows from a ucp analogue of induction. 
Let $\Delta\le\G$ and fix a section $\sigma\colon\G/\Delta\to\G$ 
so that $\G=\bigcup_{x\in\G/\Delta}\sigma(x)\Delta$. 
Let $\p\colon C^*\Delta\to\B(\hh)$ be a ucp map 
and let $F\subset\G/\Delta$ be a finite subset. 
Then, the ucp map $\p$ dilates to a ucp map 
$\Phi^F\colon C^*\G\to\B(\ell_2F\otimes\hh)$ 
given by 
$$\Phi^F(s)=\sum_{x\in F\cap s^{-1}F}
e_{sx,x}\otimes\p(\sigma(sx)^{-1}s\sigma(x))
\in\B(\ell_2F\otimes\hh)\mbox{ for }s\in\G.$$
Indeed, let $\hhh$ be the Hilbert space arising from 
$\C\G\otimes\hh$ with the inner product 
$\ip{\delta_t\otimes\eta,\,\delta_s\otimes\xi}
=(\p(s^{-1}t)\eta\mid\xi)$, 
where $\p$ is extended on $\G$ by puting $\p(s)=0$ 
for $s\in\G\setminus\Delta$. 
Clealy, $\G$ acts on $\hhh$ by left translation. 
Since $\Phi^F$ coincides with the restriction of 
this unitary representation to the subspace generated 
by $\C\sigma(F)\otimes\hh$, it is ucp. 
Now, we assume that $\G/\Delta$ is amenable and $\Delta$ has 
the property $\mathrm{(F)}$. 
Then, there is a sequence of ucp maps
$\p_n\colon C^*\Delta\to\M_{k(n)}$ such that 
$\lim\tr_{k(n)}\p_n(s)=\delta_{1,s}$ and 
$\lim\tr_{k(n)}(\p_n(s)^*\p_n(s))=1$ 
for all $s\in\Delta$. 
It is not too hard to see that the sequence $\Phi_n^{F_n}$ on $C^*\G$ 
satisfy the same conditions for a suitably chosen 
F{\o}lner sequence $F_n\subset\G/\Delta$. 

We turn to an amalgamated free product. 
Let $A_i$ ($i\in I$) be $C^*$-algebras with 
a common finite dimensional $C^*$-subalgebra $B$ 
and let $\tau_i$ be liftable traces on $A_i$ which agree on $B$. 
We have to show the free product trace $\tau$ on 
the full amalgamated free product $\ast_{i\in I}(A_i,B)$ is liftable. 
We may assume that $I=\{1,2\}$. 
Since each $\tau_i$ is liftable, there is an asymptotically 
multiplicative and asymptotically trace preserving 
ucp map $\p_i$ from $A_i$ into a full matrix algebra $D_i$. 
We claim that $\p_i$ can be chosen 
exactly trace preserving and multiplicative on $B$. 
Indeed, by amplifying the range $D_i$ if necessary, 
we may assume that there is a trace preserving $*$-homomorphism 
$\pi_i\colon B\to D_i$ such that 
$\|\p_i(u)-\pi_i(u)\|_2\approx0$ for $u\in\U(B)$ 
(cf.\ \cite{bd} for technical details). 
Note that $\U(B)$ is compact since $B$ is finite dimensional.
We define a map $\p'_i\colon A_i\to D_i$ by 
$$\p'_i(a)=\iint_{\U(B)\times\U(B)}
 \pi_i(u)\p_i(u^*av)\pi_i(v^*)\,du\,dv.$$
It is not too hard to see that $\p'_i$ is a completely positive 
$B$-bimodule map. 
By our assumption, $e_i:=1-\p'_i(1)\in \pi_i(B)'_+$ is close 
to zero in the $2$-norm. 
Finally, take a ucp extension $\tilde{\pi}_i\colon A_i\to D_i$ of $\pi_i$
and define a ucp map $\p''_i\colon A_i\to D_i$ by 
$\p''_i(a)=\p'_i(a)+e_i^{1/2}\tilde{\pi}_i(a)e_i^{1/2}$ for $a\in A_i$. 
Then, $\p''_i$ is a small perturbation of $\p_i$ and agrees 
with $\pi_i$ on $B$. 

Let $A_i^0\subset A_i$ be the orthogonal complement of $B$ w.r.t.\ 
the trace $\tau_i$.
By Boca's theorem \cite{boca}, the $B$-bimodule map 
$\Phi\colon(A_1,B)\ast(A_2,B)\to(D_1,B)\ast(D_2,B)$, 
defined by $\Phi(a_1\cdots a_n)=\p_{i_1}(a_1)\cdots\p_{i_n}(a_n)$ 
for $a_k\in A_{i_k}^0$, $i_1\neq\cdots\neq i_n$, 
is a ucp map between the full amalgamated products. 
Since each $\p_i$ is asymptotically trace preserving and 
asymptotically multiplicative, so is $\Phi$. 
Thus, we reduced the problem to the case where 
$A_i$ are full matrix algebras $D_i$. 
This case was proved by Brown and Dykema \cite{bd}. 
Indeed, under this assumption, 
they showed that the reduced amalgamated free product 
is embeddable into an interpolated free group factor 
(and a fortiori into $R^\omega$) and a $*$-homomorphism 
from the full amalgamated free product into $R^\omega$ 
is ucp liftable. 
\end{proof}

Brown suggested a possibility that all one-relator groups 
have the property $\mathrm{(F)}$. 
We note that the Baumslag-Solitar groups 
$BS(p,q)=\langle a,t \mid ta^pt^{-1}=a^q\rangle$, 
which are typical examples of non-residually finite 
one-relator groups (at least for $1<p<q$), 
are hyperlinear \cite{radulescu2}
and moreover have the property $\mathrm{(F)}$ 
since they are HNN extension of the cyclic group $\Z$ and hence 
embeddable in the semidirect products by $\Z$ 
of the increasing union of residually finite groups. 
However, it is not that all groups have the property $\mathrm{(F)}$. 
Indeed, infinite simple groups with Kazhdan's property $\mathrm{(T)}$ 
(for existence of such groups, see \cite{gromov}) 
do not have the property $\mathrm{(F)}$ 
as it was shown by Kirchberg \cite{kf}. 
We refer the reader to de la Harpe and Valette's book \cite{hv}
for Kazhdan's property $\mathrm{(T)}$.
The following proof is close to that of Kirchberg \cite{kf}, 
but uses an idea of Bekka \cite{wassermann}. 
\begin{thm}
Let $\tau$ be a trace on the full $C^*$-algebra $C^*\G$ 
of a group $\G$ with Kazhdan's property $\mathrm{(T)}$. 
Then, $\tau$ is liftable if and only if 
there is a sequence of $*$-homomorphisms 
$\pi_n\colon C^*\G\to\M_{k(n)}$ 
such that $\tau(a)=\lim_n\tr_{k(n)}\pi_n(a)$. 

In particular, a group $\G$ with the properties $\mathrm{(F)}$ 
and $\mathrm{(T)}$ is residually finite. 
\end{thm}
\begin{proof}
Let $\widehat{\G}$ be a set of irreducible representations of $\G$ 
which contains exactly one representation from each equivalence class. 
It follows that $\Theta=\bigoplus_{\theta\in\widehat{\G}}\theta$ 
is a faithful representation of $C^*\G$ on the Hilbert space 
$\hh=\bigoplus\hh_\theta$. 
We consider the faithful representation 
$(\Theta\otimes\bar{\Theta})^\infty$ of 
$C^*\G\otimes_{\min}\overline{C^*\G}$ 
on $(\hh\otimes\bar{\hh})^\infty$. 
By the assumption, the state $\mu_\tau$ 
on $C^*\G\otimes_{\min}\overline{C^*\G}$ is continuous. 
Hence, it is approximated by a vector states 
associated with a unit vector $\xi_n\in(\hh\otimes\bar{\hh})^\infty$. 
Since $\mu_\tau(s\otimes\bar{s})=1$ for all $s\in\G$, 
the sequence of unit vectors $\{\xi_n\}$ 
is almost invariant under $(\Theta\otimes\bar{\Theta})^\infty(\G)$. 
(Here, we regard $(\Theta\otimes\bar{\Theta})^\infty$ as a representation 
of $\G$.)
Hence, by Kazhdan's property $\mathrm{(T)}$, 
we may assume that $\xi_n$ are 
actually invariant under $(\Theta\otimes\bar{\Theta})^\infty(\G)$.
The well-known lemma of Schur states that for 
irreducible $\theta_1$ and $\theta_2$, 
the representation $\theta_1\otimes\bar{\theta}_2$ has 
a non-zero invariant vector if and only if 
$\theta_1$ and $\theta_2$ are equivalent and finite dimensional. 
Moreover, any invariant vector for $\theta\otimes\bar{\theta}$ 
is a constant multiple of 
$\eta_\theta:=d_\theta^{-1/2}\sum_k\zeta_k\otimes\bar{\zeta}_k$, 
where $\{\zeta_k\}_k$ is an orthonormal basis 
of the $d_\theta$-dimensional Hilbert space $\hh_\theta$. 
It follows that we may assume that the vector states associated 
with $\xi_n$ is a convex combination (with rational coefficients) 
of that of $\eta_\theta$'s with finite dimensional $\theta$. 
We obtain the conclusion by observing that 
$((\Theta(s)\otimes\bar{1})\,\eta_\theta\mid\eta_\theta)
=\tr_{d_\theta}\theta(s)$ 
for all $s\in\G$. 
\end{proof}

It is unknown whether there exists a simple property $\mathrm{(T)}$ 
group $\G$ which is hyperlinear. 
By the above theorem, the full $C^*$-algebra $C^*\G$ 
of such a group $\G$ cannot have the LLP. 
It is obvious that an inductive limit of hyperlinear groups 
is again hyperlinear. 
Since Gromov \cite{gromov} constructed infinite simple 
$\mathrm{(T)}$ groups as inductive limits 
(where connecting maps are surjective and not injective) 
of hyperbolic groups, it is particularly interesting to know 
whether all hyperbolic groups are hyperlinear 
(or more generally satisfy the property $\mathrm{(F)}$). 

\section{Further Topics and Open Problems}

Haagerup and Winsl{\o}w \cite{hw1}\cite{hw2} studied
the space $\mathrm{vN}(\hh)$ of all von Neumann algebras 
acting on a fixed separable Hilbert space $\hh$, 
equipped with the Effros-Mar{\'e}chal topology. 
Among other things, they proved that a $\mathrm{II}_1$-factor 
$M\in\mathrm{vN}(\hh)$ is embeddable into $R^\omega$ 
if and only if it is approximated by finite dimensional factors 
in $\mathrm{vN}(\hh)$. 
They also gave a new proof of equivalence between the conjectures. 

R{\u a}dulescu \cite{radulescu1}\cite{radulescu2} studied 
the space of the moments of noncommutative monomials of degree $p$ 
in a $\mathrm{II}_1$-factor $M$ and proved that it coincides 
with that for the hyperfinite $\mathrm{II}_1$-factor $R$ 
when $p\le3$. 
He also showed that the same assertion for $p=4$ is equivalent to 
the Connes' embedding problem. 

Brown \cite{brown}\cite{brown2} studied various approximation 
properties of traces on $C^*$-algebras and 
proved in particular that a separable 
$\mathrm{II}_1$-factor (or any von Neumann algebra) 
$M\subset\B(\hh)$ is embeddable into $R^\omega$ (resp.\ QWEP)
if and only if there are a weak$^*$-dense $C^*$-algbera $A$ in $M$ 
and a ucp map $\p\colon\B(\hh)\to M$ with $\p|_A=\id_A$. 
There is a related result of Brown and Dykema \cite{bd} on 
the (interpolated) free group factors $L\F_s$. 

The following questions were raised by Kirchberg \cite{klp}.
\medskip

\noindent\textbf{Problem.} 
Is the QWEP conjecture true? 
\medskip

\noindent\textbf{Problem.} 
Does there exist a non-nuclear $C^*$-algebra 
with the WEP and the LLP?
\smallskip 

What Kirchberg found is very close to an example. 
He found an short exact sequence 
$0\to\K(\ell_2)\to B\to\mathrm{cone}(C^*_{\mathrm{red}}\F_2)\to 0$ 
with $B\otimes_{\min}B^{\op}=B\otimes_{\max}B^{\op}$. 
By Theorem \ref{lance}, the $C^*$-algebra $B$ has the WEP. 
Unfortunately, it is not known whether $B$ has the LLP. 
We remark that a $C^*$-algebra with the WEP has the LLP 
if and only if it is os finitely representable in $C^*\F_\infty$ 
(cf.\ \cite{jp}).
In the remarkable paper \cite{ht}, 
Haagerup and Thorbj{\o}rnsen proved that
$C^*_{\mathrm{red}}\F_2$ is embeddable into 
$\prod\M_n/\bigoplus\M_n$, 
which has to be weakly cp complemented 
(cf.\ the remarks at the end of section 4). 
Therefore, all quasidiagonal extensions of 
$C^*_{\mathrm{red}}\F_2$ have the WEP. 
However, there are uncountably many mutually 
non-isomorphic such extensions, among which 
at most countably many can possibly have the LLP. 
(The precise statement is that the set of 
all 3-dimensional operator subspaces 
in quasidiagonal extensions of $C^*_{\mathrm{red}}\F_2$ 
is non-separable in the cb-distance topology (cf.\ \cite{jp}).)
\medskip

\noindent\textbf{Problem.} 
Does there exist a non-nuclear exact $C^*$-algebra 
with the LLP?
\smallskip

Existence of such an example is in contradiction with 
the QWEP conjecture since exact $C^*$-algebras are 
locally reflexive \cite{kcom} and locally reflexive $C^*$-algebras 
with the WEP are necessarily nuclear \cite{eh}. 
To the best of our knowledge, all $C^*$-algebras known to have 
the LLP are either nuclear or equivalent to $C^*\F_\infty$ 
in the sense that they are (cp complemented) subquotients of each other. 
We note that a $C^*$-algebra $A$ has 
both the exactness and the LLP if and only if 
$A\otimes_{\min}(\B(\ell_2)/\K(\ell_2))
=A\otimes_{\max}(\B(\ell_2)/\K(\ell_2))$. 
\medskip

\noindent\textbf{Problem.} 
Does there exist a separable $C^*$-algebra 
with the LLP, but not the LP?
\smallskip

Again, existence of such an example is in contradiction with 
the QWEP conjecture by Corollary \ref{llplp}. 
It seems that we have no example of a locally ucp liftable map 
defined on a separable $C^*$-algebra which does not have 
a global lifting.

\end{document}